%% file: main.tex
\documentclass[12pt]{article}
\usepackage{import}
\usepackage{amssymb,amsmath}
\usepackage{amsthm}
\usepackage{amsfonts,epsfig,epstopdf,titling,url,array,hyperref}
\usepackage[nameinlink]{cleveref}
\usepackage[numbers]{natbib}
\setcitestyle{comma}

\usepackage{mathrsfs}
\usepackage{mathtools}

\usepackage{tikz,graphicx,float}
\usepackage{pgfplots}
\usepgfplotslibrary{fillbetween}
\usetikzlibrary{patterns}

\pgfplotsset{compat=1.12}
\usepackage{color}

\usepackage[latin1]{inputenc}
\usepackage[T1]{fontenc}
\DeclareGraphicsExtensions{.pdf,.png,.jpg}

\newcommand{\paralambda}{\lambda}
\newcommand{\paragamma}{\gamma}

\newcommand{\bdplambda}{\beta}

\newcommand{\expect}[3][]{\mathbb{E}^{#1}_{#2}\left[#3\right]}
\newcommand{\prob}[3][]{\mathbb{P}_{#1}^{#2}\left[#3\right]}
\newcommand{\N}{\mathbb{N}} 
\newcommand{\R}{\mathbb{R}} 
\newcommand{\smallo}[1]{o\left(#1\right)}

\newcommand{\polyorder}[1]{\langle #1 \rangle}
\newcommand{\abs}[1]{\left\vert #1 \right\vert}
\newcommand{\labs}{\Big\rvert}

\newcommand{\toggle}[1]{}

\renewcommand{\exp}[1]{\textnormal{exp}\left\{#1\right\}}
\newcommand{\product}[1]{\left\langle #1\right\rangle}
\newcommand{\twobytwofig}[5]{
    \begin{figure}[H]
\begin{minipage}[t]{0.33\linewidth}
\centering
\includegraphics[scale=1, width=0.93\linewidth]{#2}
\end{minipage}%
\begin{minipage}[t]{0.33\linewidth}
\centering
\includegraphics[scale=1, width=0.9\linewidth]{#3}
\end{minipage}%
\begin{minipage}[t]{0.33\linewidth}
\includegraphics[scale=1, width=0.9\linewidth]{#4}
\end{minipage}
\caption{#1}
\label{#5}
\end{figure}
}
\newcommand{\bigthreefig}[5]{
    \begin{figure}[H]
\includegraphics[scale=0.58]{#2}
\includegraphics[scale=0.57]{#3}\\
\centering
\includegraphics[scale=0.6]{#4}\\
\caption{#1}
\label{#5}
\end{figure}
}
\newtheorem{thm}{Theorem}
\newtheorem{lem}[thm]{Lemma}
\newtheorem{prop}[thm]{Proposition}

\numberwithin{thm}{section} \numberwithin{equation}{section}

\begin{document}
\title{
Extinction time of stochastic SIRS models with small initial size of the infected population
}
\date{}
\author {Jingran Zhai\thanks{Queen Mary, University of London, email: j.zhai@qmul.ac.uk}} \maketitle

\begin{abstract}
\noindent 
The stochastic SIRS model is a continuous-time Markov chain modelling the spread of infectious diseases with temporary immunity, in a homogeneously-mixing population of fixed size $N$. We study the scaling behaviour of the extinction time of stochastic SIRS models as $N$ tends to infinity. When the initial size of infected population is small, we obtain the closed-form expression of the asymptotic distribution of this extinction time, and compare it with the data from Monte Carlo simulation.
\end{abstract}

\noindent
MSC:
\vskip0.2cm
\noindent
{\it Keywords:} 
\noindent

\noindent
\section{Introduction}
Stochastic epidemic models play important roles in understanding the dynamic of infectious diseases. The criticality of a stochastic epidemic model is a threshold in transmission rate such that a subcritical epidemic tends to die out quickly and a supercritical epidemic tends to prevail in the population. Over the last few years, the behaviours of near-critical epidemics have drawn a lot of attention. As \citet{britton2015five} pointed out, many diseases, especially those under eradication campaign, are near-critical, and therefore understanding such behaviour is a significant challenge in stochastic epidemic modelling. Analysis on the stochastic SIS \citep{Nasell1996, dolgoarshinnykh2006critical, 2018Foxall} and stochastic SIR model \citep{dolgoarshinnykh2006critical} has shown that within certain near-critical regime, the epidemics described by these models demonstrate persistence behaviours different from strongly sub- or supercritical cases. 

We aim to extend this analysis to the stochastic SIRS model. In particular, we investigate the extinction time, which is the time for a population with given initial state to reach a state of `zero cases'. 

Estimating the extinction time of stochastic epidemic models as a problem itself has drawn a lot of attention. \citet{barbour1975duration} obtained the asymptotic distribution of the stochastic SIR model. Chronologically, \citet{KryscioLefvre1989,anderson1992infectious, Nasell1996,Doering2005} all studied the expectation of the SIS extinction time. Recently, the asymptotic distribution of the SIS extinction time is obtained by \citet{Malwina2013} and \citet{2018Foxall}. As far as we are aware, the only available result regarding the stochastic SIRS model is by \cite{Dolgoarshinnykh2003}, who obtained the expected extinction time for strongly supercritical SIRS model with initial size of infected and immune population both being of order $N$ (the size of the total population).

From a practical point of view, many human infections have a temporary but significant duration of immunity and thus are better modelled by an SIRS model, especially when the observation time window is long. From a mathematical point of view, the SIRS model introduces extra complexity comparing to the SIS and SIR models, and thus the extension is non-trivial: while the SIS and the SIR models are both driven by a single parameter (the transmission rate), the SIRS model incorporates a second parameter describing the average duration of immunity. The SIS model is simpler since it is a one-dimensional process and its mean-field ODE approximation has an explicit solution. Despite being two-dimensional, the SIR model has a monotonicity which simplifies the analysis. In addition, there is no explicit solution to the ODE system describing deterministic SIRS models.

As the first step of solving this problem, in this paper we focus on the behaviour of the extinction time when the initial size of infected population is small. This is usually the case when the infection is introduced to a new population, or when the epidemic is approaching extinction under intervention. The precise definition of `small' is complicated and depends on the parameters, which will be made clear in our main result \Cref{thm:small_ini_main_result}. The following is an example of the scenario we shall investigate: the strongly subcritical SIRS model ($\paralambda<1$ and $\paragamma>0$) with $I_0 = R_0 = N^{1/2}$. We are able to obtain the closed-form expression of the asymptotic distribution of the extinction time with various near-critical `parameter and initial state' combinations. \\

The stochastic SIRS model describes the spread of a disease with no incubation period and a temporary immunity in a closed population of size $N$. A population is said to be \textit{closed} if it has no birth, death, immigration or migration. We assume that each infected individual contacts any other individual at rate $\paralambda/N$ and will transmit the disease if his/her contact is susceptible. Parameter $\paralambda\in\R_+$ is known as the \textit{transmission rate}. Once infected, each individual is immediately infectious and will recover at rate $1$ independently. Each recovered individual loses immunity at rate $\paragamma\in\R_+$ and becomes susceptible independently. To study the near-critical parameter regime, we assume that $\paralambda$ and $\paragamma$ are bounded and depend on $N$. We say that the stochastic SIRS model is strongly subcritical if $\paralambda<1$ as $N\to \infty$ and is near-critical if $\paralambda\to 1$ as $N\to \infty$. This is consistent with the criticality of the stochastic SIS and SIR models.

Formally, the stochastic SIRS model is defined as a two-dimensional continuous-time Markov chain $(I^N,R^N)$ valued in $\{(i,r)\in\N^2: 0 \leq i+r\leq N\}$, with the following transition rates:
\begin{alignat}{2}
    &(i,r) \to (i+1,r),\, && \textnormal{ at rate } \paralambda (N-i-r)i/N,\nonumber\\
    &(i,r) \to (i,r-1),\, && \textnormal{ at rate } \paragamma r,\label{def:SIRS}\\
    &(i,r) \to (i-1,r+1),\, && \textnormal{ at rate } i.\nonumber
\end{alignat}

It has been long noticed that the trajectory of the size of infected population $I^N$ can be well-approximated by linear birth-death processes when $I^N$ is small. Such approximation is done by constructing an order-preserving coupling between birth-death processes. Among the existing works that use this technique, \citet{barbour1975duration}, \citet{Malwina2013} and \citet{2018Foxall} all described a version of the construction of this coupling. The various constructions are the same in nature, as described in \Cref{appd:coupling_for_bdp}.  \citet{barbour1975duration} studied the stochastic SIR model with transmission rate $\paralambda$ independent of $N$, and chronologically \citet{Malwina2013} and \citet{2018Foxall} both used this technique on the subcritical and near-critical stochastic SIS model.  In particular, the discussion made by \citet{2018Foxall} is the most comprehensive, in the sense that it covers all possible scenarios of the stochastic SIS model where this technique is applicable. Our work is motivated by \citet{2018Foxall}.

\section{Main results}

Throughout this paper, we use the following asymptotic notations:\\
for functions $f(x)$ and $g(x)$, 
\begin{itemize}
    \item if there exists constant $K>0$ s.t. 
$K^{-1}|g(x)|\leq |f(x)|\leq K|g(x)|,$
then we say $f(x)\asymp g(x)$;
    \item if 
$\lim_{x\to\infty}\frac{f(x)}{g(x)}= 0$,
then we say $f(x)=o(g(x))$, $f(x)\ll g(x)$ or $g(x)\gg f(x)$;
    \item if $\lim_{x\to\infty}\frac{f(x)}{g(x)}= 1$,
then we say $f(x)\sim g(x)$.
\end{itemize}

For each $N$, the extinction time of the stochastic SIRS model $(I^N,R^N)$ is defined as $T^N_o:=\inf\{t:I^N_t = 0\}$. 

\begin{thm}\label{thm:small_ini_main_result}
Consider a sequence of stochastic SIRS models $(I^N,R^N)$ indexed by $N\in\N$, with parameters $\paralambda = \paralambda(N)>0$ and $\paragamma=\paragamma(N)> 0$, and initial states $(I^N_0,R^N_0) = (I_0(N),R_0(N))$. 

If one of the following conditions is satisfied, then we have the closed-form expression of the asymptotic distribution of $T^N_o$. 

Cases 1.1-1.3 are cases where both the initial size of infected population $I_0$ and immune population $R_0$ are small.
\begin{itemize}
    \item[] \textbf{Case 1.1: }$I_0|1-\paralambda|\to 0$, $I_0R_0 = o(N)$, $I_0 = o(N^{1/2}\paragamma^{1/2})$.
    
    If $I_0 = O(1)$,
    \[\lim_{N\to\infty}\prob{N}{T^N_{o}\leq w} = \left(1+\frac{1}{w}\right)^{-I_0};\]
    and if $I_0 \to\infty$,
    \[\lim_{N\to\infty}\prob{N}{\frac{T^N_{o}}{I_0}\leq w} = e^{-\frac{1}{w}}.\]
    \item[] \textbf{Case 1.2: }$I_0(1-\paralambda)\to a >0$, $\paralambda=\paralambda(N) < 1$, and $I_0 = \smallo{N^{1/2}\paragamma^{1/2}}$, $I_0R_0 = \smallo{N}$. 
    
    If $I_0 = O(1)$, 
    \[\lim_{N\to\infty}\prob{N}{T^N_{o}\leq w} = \left(1+\frac{a}{e^{a w}-1}\right)^{-I_0};\]
    and if $I_0 \to\infty$,
     \[\lim_{N\to\infty}\prob{N}{\frac{T^N_{o}}{I_0}\leq w} = \exp{-\frac{a}{e^{aw}-1}}.\]
    \item[] \textbf{Case 1.3: }$I_0(1-\paralambda)\to \infty$, $\paralambda=\paralambda(N) < 1$, $I_0 = \smallo{\frac{N(1-\paralambda)\paragamma}{\log I_0(1-\paralambda)}}$, and $R_0\log I_0(1-\paralambda) = \smallo{N(1-\paralambda)}$. Then
    \[\lim_{N\to\infty}\prob{N}{(1-\paralambda)T^N_{o}- \log (1-\paralambda)I_0\leq w} = e^{-e^{-w}}.\]
    \end{itemize}
Cases 2.1 and 2.2 are cases where $I_0$ is small and $R_0$ is of order $N$. 
\begin{itemize}
    \item[] \textbf{Case 2.1: } $I_0 = O(1)$, $R_0 = r_0N$, $r_0>0$, $\paralambda=\paralambda(N)\leq 1$ and $\paragamma = o(1)$.
Let $a:=\lim_{N\to\infty}\left(1-\paralambda+\paralambda r_0\right)$, then
    \begin{align*}
        \lim_{N\to\infty}\prob{N}{T^N_o\leq  w}   = \left(1+\frac{a}{e^{aw}-1}\right)^{-I_0}.
    \end{align*}
    \item[] \textbf{Case 2.2: } $I_0\to\infty$, $R_0 = r_0N$, $r_0>0$, $\paralambda=\paralambda(N)\leq 1$, and there exist $\epsilon_1,\epsilon_2>0$ such that $I_0 = \smallo{N^{1-\epsilon_1}}$ and $\paragamma = \smallo{N^{-\epsilon_2}}$. Let $a:=\lim_{N\to\infty}\left(1-\paralambda+\paralambda r_0\right)$, then
    \begin{align*}
        \lim_{N\to\infty}\prob{N}{aT^N_o -\log(aI_0)\leq w} & = e^{-e^{-w}}.
    \end{align*}
\end{itemize}
\end{thm}
The cases above cover completely the parameter regime $\{(\paralambda,\paragamma):\paralambda\leq 1\}$, and Case 1.1 also covers a subset of the parameter regime $\{(\paralambda,\paragamma):\paralambda\geq 1\}$. Given $(\paralambda,\paragamma)$, we illustrate the conditions with respect to the initial states of each case in \Cref{thm:small_ini_main_result} using the following diagram.

We use $\polyorder{f(N)}$ to denote the scaling of $f(N)$: for any function $f = f(N)$, $\polyorder{f} = a \in\R$ if and only if $\abs{f(N)}\asymp N^a$. If $\abs{f(N)}$ tends to infinity faster then any polynomials, we say $\polyorder{f} = \infty$; and if $\abs{f(N)}$ tends to infinity slower then any polynomials, we say $\polyorder{f} = 0+$ and  $\polyorder{1/f} = 0-$.

In \Cref{fig:graph_sub_small}, the first diagram illustrates the location of the initial states of all five possible cases in \Cref{thm:small_ini_main_result} when
\[\bigg\{(\paralambda,\paragamma):N^{1/2}(1-\paralambda)\paragamma^{1/2}\to \infty\bigg\},\]  
and the second diagram illustrates the location of the initial states of the three possible cases when
\[\bigg\{(\paralambda,\paragamma):N^{1/2}(1-\paralambda)\paragamma^{1/2}<\infty\bigg\}.\] 

\begin{figure}[H]
  \centering
  \scalebox{.8}{
    \import{images/}{graph_sub.tex}
    }
    \caption{Diagram of the initial state conditions given parameters $(\paralambda, \paragamma)$, where $a = - \polyorder{1-\paralambda}$, $b = \frac{1+\polyorder{\paragamma}}{2}$, $c = 1+\polyorder{1-\paralambda}+\polyorder{\paragamma}$. The numbers denote the cases in \Cref{thm:small_ini_main_result}.}
     \label{fig:graph_sub_small}
\end{figure}

\section{Properties of linear birth-death-immigration processes}
By \textit{linear birth-death-immigration process}, we mean the continuous-time Markov chains $(L_t)_{t\geq 0}$ valued in $\N$, with transition rates 
\begin{alignat}{2}
    &x \to x+1,\, && \textnormal{ at rate } \bdplambda x+\alpha,\nonumber\\
    &x \to x-1,\, && \textnormal{ at rate } \mu x,\nonumber
\end{alignat}
where parameters $\bdplambda, \mu\geq 0$ are known as \textit{birth rate} and \textit{death rate} respectively, and $\alpha\geq 0$ is known as \textit{immigration rate}. To prove \Cref{thm:small_ini_main_result}, we need some properties and estimations of linear birth-death processes ($\alpha =0$) and linear immigration-death processes ($\bdplambda = 0$).

The first result is the closed-form expression of the limit distribution of the extinction time of linear birth-death processes. The distribution itself is well-known (See e.g. \citet{renshaw2015stochastic}) but its limit form when the parameters and the initial states both depend on $N$ was obtained recently by \citet{2018Foxall}.
\begin{thm}[Theorem 1,  \citep{2018Foxall}]\label{thm:bdp_extinction_limit}
Let $\{L^N\}_{N\in\N}$ be a sequence of linear birth-death processes with birth rate $\bdplambda_N>0$ and death rate $\mu_N > 0$.

Let $T^N_{bdp} := \inf\left\lbrace t\geq 0: L^N_t = 0\right\rbrace$. The distribution of $T^N_{bdp}$ converges to the following limits as $N\to\infty$:\\
Suppose that $L^N_0 = L_0$ is a constant independent of $N$:
\begin{enumerate}
    \item If $\mu_N-\bdplambda_N\to 0$, then 
\[\lim_{N\to\infty}\prob[L_0]{N}{T^N_{bdp}\leq w} = \left(1+\frac{1}{w}\right)^{-L_0},\quad w>0;\]
    \item If $\mu_N-\bdplambda_N\to a>0$, then
\[\lim_{N\to\infty}\prob[L_0]{N}{T^N_{bdp}\leq w} = \left(1+\frac{a}{e^{a w}-1}\right)^{-L_0},\quad w>0.\]
\end{enumerate}
Suppose that $L^N_0 = L_0(N)\to \infty$:
\begin{enumerate}
\setcounter{enumi}{2}
    \item If $L_0(\mu_N-\bdplambda_N)\to 0$, then 
\[\lim_{N\to\infty}\prob[L_0]{N}{\frac{T^N_{bdp}}{L_0}\leq w} = e^{-\frac{1}{w}},\quad w>0;\]
    \item If $L_0(\mu_N-\bdplambda_N)\to a>0$, then 
\[\lim_{N\to\infty}\prob[L_0]{N}{\frac{T^N_{bdp}}{L_0}\leq w} = \exp{-\frac{a}{e^{aw}-1}},\quad w>0;\]
    \item If $L_0(\mu_N-\bdplambda_N)\to \infty$, then 
\[\lim_{N\to\infty}\prob[L_0]{N}{(\mu_N-\bdplambda_N)T^N_{bdp}- \log L_0(1-\bdplambda_N/\mu_{N})\leq w} = e^{-e^{-w}},\quad w\in\R.\]
In particular, if for some $a(N)\sim L_0(\mu_N-\bdplambda_N)$, we have $L_0(\mu_N-\bdplambda_N)-a(N) = \smallo{\frac{a}{\log a}}$, then we can also write the limit distribution as
\[\lim_{N\to\infty}\prob[L_0]{N}{a(N)\frac{T^N_{bdp}}{L_0}- \log \frac{a(N)}{\mu_{N}}\leq w} = e^{-e^{-w}}.\]
    \end{enumerate}
\end{thm}

Next, we need some auxiliary results to estimate the probability of some birth-death-immigration processes hitting s certain barrier. Although the approach is routine, since we assume in addition that the parameters and the initial states are of various scaling of $N$, we will state the full proofs below. 
\begin{lem}[Hitting probability of linear birth-death processes]\label{lem:hitting_lbd}
Let $L = (L_t)_{t\geq 0}$ be a linear birth-death process with birth rate $\bdplambda=\bdplambda(N)>0$ and death rate $1$, and $L_0 = l(N)\in\N$. 

As $N\to\infty$, the probability that $L_t$ ever reaches $k(N)>l(N)$ tends to $0$, if one of the following conditions holds:
\begin{enumerate}
    \item $(1-\bdplambda)( k-l)\to \infty$, $\bdplambda(N)<1$ for all $N\in\N$ and $\lim_{N\to\infty}\bdplambda(N)\leq 1$;
    \item $l(N)= o(k(N))$, $(1-\bdplambda)k \to 0$ and $\bdplambda\to 1$.
\end{enumerate}
\end{lem}
\begin{proof}
Let $h_i$ be the probability of $L$ ever hitting $k$ from $L_0 = i$, $i\in\N$. Then $\{h_i\}_{i\leq k}$ is the minimal non-negative solution (See Theorem 3.3.1, p.112, \cite{norris1998markov}) of 
\begin{align*}
    0 & = \bdplambda (h_{i+1}-h_i) + (h_{i-1}-h_i),\quad 1<i<k,\\
    h_k & = 1,\quad h_0=0.
\end{align*}
It has a unique solution 
\[ h_i=\frac{\bdplambda^{-i}-1}{\bdplambda^{-k}-1},\quad i\leq k.\]

If $(1-\bdplambda)( k-l)\to \infty$, $\bdplambda(N)<1$ for all $N\in\N$, and $\lim_{N\to\infty}\bdplambda(N)\leq 1$, we have
\[h_l\leq \bdplambda^{k-l} = \left((1-(1-\bdplambda))^{\frac{1}{1-\bdplambda}}\right)^{(1-\bdplambda)(k-l)}\to 0.\]

If $(1-\bdplambda)k \to 0$ and $\bdplambda\to 1$, we have
\[\bdplambda^{-i} = \left((1-(1-\bdplambda))^{\frac{1}{1-\bdplambda}}\right)^{(1-\bdplambda)i}\to 1, \quad i\leq k.\]
Notice that this is true even if $\bdplambda(N)\geq 1$ for some $N\in\N$.

It follows that
\[\lim_{N\to\infty} h_l = \lim_{N\to\infty}\frac{\sum_{i=0}^{l-1}\bdplambda^{-i}}{\sum_{i=0}^{k-1}\bdplambda^{-i}} = \lim_{N\to\infty}\frac{l-1}{k-1}= 0,\]
when $l = o(k)$. 
\end{proof}
\begin{lem}[Hitting probability of immigration-death processes absorbing at $0$]\label{lem:hitting_id}
For any given $N\in\N$, let $L= (L_t)_{t\geq 0}$ be an immigration-death process with immigration rate $\alpha(N)>0$ and death rate $\mu(N)>0$, absorbing at $0$. That is, $L$ has the transition rates for $x\geq 1$ as follows:
    \begin{alignat}{2}
        &x \to x+1,\, && \textnormal{ at rate }  \alpha,\nonumber\\
        &x \to x-1,\, && \textnormal{ at rate } \mu x;\nonumber
    \end{alignat}
and $L$ remains at $0$ once it hits $0$. 

Let $L_0 = l(N)\to\infty$. 
Then the probability that $L_t$ ever reaches $2l(N)$ tends to $0$ as $N\to\infty$, if 
\[\lim_{N\to\infty}\frac{l(N)\mu(N)}{\alpha(N)}> e.\]
\end{lem}
\begin{proof}
    Let $h_i$ be the probability of $L$ ever hitting $k$ from $L_0 = i$, $i\leq 2l$. Then $\{h_i\}_{i\leq 2l}$ is the minimal non-negative solution of 
\begin{align*}
    0   & = \alpha (h_{i+1}-h_i) + \mu i(h_{i-1}-h_i),\quad 1\leq i<2l,\\
    h_{2l} & = 1,\quad h_0 = 0.
\end{align*}
It has a unique solution $\{x_i\}_{i\leq 2l}$, where
\begin{align*}
    x_i & =  \frac{\sum_{k=0}^{i-1}(\mu/\alpha)^{k}k!}{\sum_{k=0}^{2l-1}(\mu/\alpha)^{k}k!}+h_0\left(1- \frac{\sum_{k=0}^{i-1}(\mu/\alpha)^{k}k!}{\sum_{k=0}^{2l-1}(\mu/\alpha)^{k}k!} \right),\quad 1\leq i\leq 2l,\\
    x_0 &= h_0=0,
\end{align*}
and thus
\begin{align*}
    \frac{\sum_{k=0}^{l-1}(\mu/\alpha)^{k}k!}{\sum_{k=0}^{2l-1}(\mu/\alpha)^{k}k!} =  h_l.
\end{align*}
When $\lim_{N\to\infty}\frac{l\mu }{e\alpha}>1$, we have 
\[\sum_{k=l}^{2l-1}(\mu/\alpha)^{k}k!\geq (\mu/\alpha)^{l}l!\sum_{k=0}^{l-1}(\mu/\alpha)^{k}k!,\]
and for sufficiently large $N$,
\begin{align*}
    h_l \leq \frac{\sum_{k=0}^{l-1}(\mu/\alpha)^{k}k!}{\sum_{k=0}^{l-1}(\mu/\alpha)^{k}k!+(\mu/\alpha)^{l}l!\sum_{k=0}^{l-1}(\mu/\alpha)^{k}k!}
    =  \frac{1}{1+(\mu/\alpha)^{l}l!} 
    \leq \frac{1}{1+l^{1/2}\left(\frac{l\mu }{e\alpha}\right)^{l}\sqrt{2\pi}} \to 0,
\end{align*}
where we use a well-known bound for the factorial function $l!> \sqrt{2\pi}l^{l+1/2}e^{-l}$ from \citet{stirlingformula}.
\end{proof}
\begin{lem}[Hitting probability of immigration-death processes]\label{lem:hitting_id_nonabsorb}
For given $N\in\N$, let $L= (L_t)_{t\geq 0}$ be an immigration-death process with immigration rate $\alpha(N)>0$ and death rate $\mu(N)>0$. That is, $L$ has the transition rates for $x\geq 0$ as follows:
    \begin{alignat}{2}
        &x \to x+1,\, && \textnormal{ at rate }  \alpha,\nonumber\\
        &x \to x-1,\, && \textnormal{ at rate } \mu x.\nonumber
    \end{alignat}
If $L_0 = l(N)\to\infty$, $\mu = O(1)$ and $\alpha = o(l\mu)$, then for $t_0 = t_0(N)\to \infty$ satisfying $t_0 =\smallo{ \left(\frac{l\mu}{\alpha e}\right)^{l}}$, the probability of the event `$L_t$ reaches $2l$ before $t=t_0$' tends to $0$ as $N\to\infty$. 
\end{lem}
\begin{proof}
Notice that $L_0 = l(N)\to\infty$, $\mu = O(1)$ and $\alpha = o(l\mu)$ imply 
\[\lim_{N\to\infty}\frac{l(N)\mu(N)}{\alpha(N)}> e.\]
Under this condition, in \Cref{lem:hitting_id}, we have estimated the probability for $L_t$ starting from $l$ to ever reach $2l$ before reaching $0$, denoted as $h_l$, and have 
\[h_l = \smallo{\left(\frac{l\mu}{\alpha e}\right)^{-l}}.\]
To prove the statement in \Cref{lem:hitting_id_nonabsorb}, we argue that with probability tending to $1$, $L_t$ can reach $l$ from $0$ at most $\lceil t_0\rceil $ times within time interval $[0,t_0]$. If this is indeed the case, then 
\[\prob[]{}{\sup_{t\in[0,t_0]}L_t\geq 2l}\leq (\lceil t_0\rceil+1) h_l\to 0.\]

Since $L_t$ is stochastically dominated by Poisson process $C_t$ with rate $\alpha$, by order-preserving coupling, we have
\[\prob[]{}{L_t \text{ travels from }0 \text{ to }l \text{ at least}\lceil t_0\rceil \text{ times within }[0,t_0]}\leq \prob[]{}{C_{t_0}\geq l\lceil t_0\rceil },\]
where $C_{t_0}\sim \text{Poisson}(\alpha t_0)$.
Since $l\lceil t_0\rceil > \alpha t_0$ for all sufficiently large $N$, we have the following bound on the tail probability of Poisson distributions (Theorem 5.4, p.97, \cite{mitzenmacher2005probability})
\[\prob[]{}{C_{t_0}\geq l\lceil t_0\rceil}\leq \left(\frac{e\alpha t_0}{l\lceil t_0\rceil}\right)^{l\lceil t_0\rceil}e^{-\alpha t_0}\to 0.\]
\end{proof}
The next result concerns finding an upper bound for $\int_0^t I^N_sds$ for $t>0$, where $I^N$ is the first component of a stochastic SIRS model. By order-preserve coupling, it suffices to estimate the same integral of a linear birth-death process which dominates $I^N$.
\begin{lem}\label{lem:integral_BDP}
Let $L(l)=(L_t)_{t\geq 0}$ be a linear birth-death process with birth rate $\bdplambda =\bdplambda(N)>0$, death rate $\mu = \mu(N)>\bdplambda(N)$, and $L_0=l(N)$ for $N\in\N$.

Let
\[T_L :=\inf\left\lbrace t: L_t = 0\right\rbrace,\]
and
\[H(l) := \int_0^{T_L}L_sds.\]
Then for each $N\in\N$, $L_0 = l(N)$ and $\delta = \delta(N)>0$, we have
\begin{align}\label{eq:integralBDP_bound}
\prob[]{}{H(l)>\delta} \leq \frac{l}{(\mu-\bdplambda)\delta}.
\end{align}
\end{lem}
\begin{proof}
We fix $N$ throughout the proof.

Denote the Laplace transform of $H$ with $L_0 = l(N)$ as 
\[H^*(a;l) := \expect{}{e^{-a H(l)}},\quad a\geq 0.\]
Let $S$ denote the sojourn time of $L(l)$ before its first jump. The explicit expression of $H^*(a;l)$ can be obtained following a first-step analysis (e.g. p.482, \cite{mcneil1970integral}):
\begin{align*}
H^*(a;l)& = \frac{\bdplambda}{\bdplambda+\mu} \expect{}{e^{-a (H(l+1)+Sl)}}+\frac{\mu}{\bdplambda+\mu} \expect{}{e^{-a (H(l-1)+Sl)}},\\
 & = \frac{\bdplambda}{\bdplambda+\mu}H^*(a;l+1)\int^{\infty}_0e^{-l a s}(\bdplambda+\mu)le^{-(\bdplambda+\mu)ls}ds \\
& + \frac{\mu}{\bdplambda+\mu}H^*(a;l-1)\int^{\infty}_0e^{-l a s}(\bdplambda+\mu)le^{-(\bdplambda+\mu)ls}ds,
\end{align*}
Then
\begin{align*}
H^*(a;l) = \left(\bdplambda H^*(a;l+1)+\mu H^*(a;l-1)\right)(\bdplambda+\mu+ a)^{-1}.
\end{align*}
The solution of the above is
\[H^*(a;l) = \left(\frac{\bdplambda+\mu+a - \sqrt[]{(\bdplambda+\mu+a)^2-4\bdplambda\mu }}{2\bdplambda}\right)^l, \quad l\geq 1.\]
\[\expect{}{H(l)} = -\frac{dH^*(a;l)}{da}\Big\vert_{a = 0} = \frac{l}{\mu-\bdplambda}.\]
By the Markov inequality, we have for each $N\in\N$ and any $\delta = \delta(N)>0$, 
\begin{align*}
\prob[]{}{H(l)>\delta}\leq \frac{\expect{}{H(l)}}{\delta} = \frac{l}{(\mu-\bdplambda)\delta}.
\end{align*}
\end{proof}

\section{Proof of the main results}
Now we are ready to prove our main result \Cref{thm:small_ini_main_result}.

The process $I^N$ has the following transition rates at time $t$ when $I^N_t = x$:
\begin{alignat}{2}
    &x \to x+1,\, && \textnormal{ at rate } \paralambda\left(1-N^{-1}(x+R^N_t)\right) x,\nonumber\\
    &x \to x-1,\, && \textnormal{ at rate } x.\nonumber
\end{alignat}
 The general idea of the proof is that we will sandwich $I^N$ between two linear birth-death processes whose extinction times have the same asymptotic distributions, according to \Cref{thm:bdp_extinction_limit}. The construction of such coupling follows from \Cref{appd:coupling_for_bdp}. To make sure the birth rates and death rates are of the correct order, we will need to find upper-bounds held with high probability for $I^N$ and $R^N$.\\

The intuition behind discussing two broad scenarios depending on the order of $R_0$ is as follows: 

If $R_0(N)/N\to 0$, then $I^N$, with small initial value and additional assumptions, will have a birth rate close to $\paralambda$. Looking at \Cref{thm:bdp_extinction_limit}, it makes sense to discuss three different cases within this scenario based on the limit of $I_0(1-\paralambda)$. 

If $R_0(N)/N\to r_0 \in (0,1]$, then $I^N$, with small initial value and additional assumptions, will have a birth rate close to $\paralambda(1-r_0)$. Depending on whether $I_0 = O(1)$, we can divide this scenario into two cases corresponding to the last two cases in \Cref{thm:bdp_extinction_limit}. \\

\subsection{Proof of Cases 1.1-1.3}
The proof of Cases 1.1 to 1.3 follows the same idea: we choose an appropriate $k(N)$ and $m(N)$ such that, $R_0\leq m(N)$ for sufficiently large $N$, and as $N\to\infty$,
\[\prob{}{R^N_t\leq 2m(N), I^N_t\leq k(N), \,\forall t\geq 0}\to 1.\]
Define two linear birth-death processes $\underline{L}$, $\overline{L}$, such that $\overline{L}$ has birth rate $\paralambda$ and death rate $1$, and $\underline{L}$ has birth rate $\paralambda(1-\frac{k(N)+2m(N)}{N})$ and death rate 1. Let $\underline{L}_0=\overline{L}_0 = I_0$. Then we only need to check that the extinction times $T_{\underline{L}}$ and $T_{\overline{L}}$ have the same asymptotic distributions.

We will state the proof of Case 1.1 in full details, and omit the repeated content in Cases 1.2 and 1.3.\\ 

\textbf{Case 1.1: }$I_0|1-\paralambda|\to 0$, $I_0R_0 = o(N)$, $I_0 = o(N^{1/2}\paragamma^{1/2})$.

Notice in this case it is necessary that $|1-\paralambda| \to 0$.

Since $I_0 = o(N^{1/2}\paragamma^{1/2})$, we can find $\kappa(N)\to\infty$ such that 
\[\kappa(N)\ll  \left(N^{1/2}\paragamma^{1/2}\wedge \abs{1-\paralambda}^{-1}\right)I_0^{-1} .\]
Let $k(N) := I_0\kappa(N)$. 
There is an order-preserving coupling between $I^N$ and $\overline{L}$ such that $I^N_t\leq \overline{L}_t$, for all $t\geq 0$.
Since $I_0\ll k$, $(1-\paralambda)k(N)\to 0$ and $\paralambda\to 1$, we can apply the second case in \Cref{lem:hitting_lbd} to $\overline{L}$, and obtain that with probability tending to $1$, $I^N_t\leq \overline{L}_t\leq k(N)$.\\ 

For $N\in\N$, conditioned on $\{I^N_t\leq k(N)\}$, each $R^N$ is stochastically dominated by an immigration-death process $M=(M_t)_{t\geq 0}$ with immigration rate $k(N)$ and death rate $\paragamma$ and $M_0 \geq  R_0$. \\

Let $M_0 = m(N) := N^{1/2}\paragamma^{-1/2}\vee R_0$, and $t_0 = M_0\paragamma$. It is obvious that $M_0\to \infty$ and $k= o(M_0\paragamma)$.  

Since $\frac{M_0\paragamma}{ke}\to\infty $ and $M_0\paragamma=O(M_0)$, we have
\[t_0=\smallo{\left(\frac{M_0\paragamma}{ke}\right)^{M_0}}.\]

Thus all the conditions of \Cref{lem:hitting_id_nonabsorb} are met, and we have
\[\prob{}{R^N_t\geq 2m(N), \,\forall t\leq t_0 \, \labs \, I^N_t\leq k(N)}\leq \prob{}{M_t\geq 2M_0, \,\forall t\leq t_0} \to 0.\]
It follows that with probability tending to $1$,
\[\paralambda\geq \paralambda\left(1-\frac{I^N_t+R^N_t}{N}\right)\geq \paralambda\left(1-\frac{k(N)+2m(N)}{N}\right).\]

Denote $T_{\underline{L}}:=\inf\{t:\underline{L}_t = 0\}$ and $T_{\overline{L}}:=\inf\{t:\overline{L}_t = 0\}$.
From Cases 1 and 3 of \Cref{thm:bdp_extinction_limit}, we have $T_{\overline{L}}$ is of order $I_0 = o(t_0)$. It follows that as $N\to\infty$, 
\[ \prob{}{T_{\overline{L}} < t_0}\to 1.\]

For each $N\in\N$, conditioned on 
\[\bigg\{I^N_t\leq 2I_0,\, R^N_t\leq 2m(N), \,\forall t\leq t_0 \bigg\},\]
there is an order-preserving coupling between $\underline{L}$ and $I^N$ and between $I^N$ and $\overline{L}$ such that $\underline{L}_t\leq I^N_t \leq \overline{L}_t$ for all $t\geq 0$. For sufficiently large $N$, we have
\[\prob{}{T_{\underline{L}}\leq T^N_o \leq T_{\overline{L}}<t_0}\geq \prob{}{I^N_t\leq k(N), R^N_t\leq 2m(N), \,\forall t\leq t_0}.\]

Notice that $I_0M_0 \leq I_0N^{1/2}\paragamma^{-1/2}+I_0R_0 =  o(N)$. Since 
\begin{align*}
     \lim_{N\to\infty}\left(1-\paralambda\left(1-\frac{(k(N)+2m(N))}{N}\right)\right)\underline{L}_0 & = \lim_{N\to\infty}(1-\paralambda)I_0+ \lim_{N\to\infty}\paralambda\frac{I_0(k(N)+2m(N))}{N}\\
     & =\lim_{N\to\infty}(1-\paralambda)\overline{L}_0 = 0,
\end{align*}
the asymptotic distribution of $T_o^N$ follows from Case 1 in \Cref{thm:bdp_extinction_limit} if $I_0 = O(1)$, and Case 3 if $I_0\to \infty$. \\

\textbf{Case 1.2: }$I_0(1-\paralambda)\to a>0$ and $I_0 = \smallo{N^{1/2}\paragamma^{1/2}}$, $I_0R_0 = \smallo{N}$. 

Notice that this is only possible if 
\begin{align}\label{eq:para_sub1}
    (1-\paralambda)N^{1/2}\paragamma^{1/2}\to\infty.
\end{align}
This case covers the scenarios where $\paralambda$ is independent of $N$ and $I_0=O(1)$.

Let $m(N) := N^{1/2}\paragamma^{-1/2}\vee R_0\to\infty $. Since $I_0 =\smallo{ m \paragamma}$, by letting $k(N) := \sqrt{I_0 m\paragamma}\to\infty$, we have \[(1-\paralambda)\sqrt{I_0m\paragamma}\gg (1-\paralambda)I_0 \asymp 1 .\] 
By the first case in \Cref{lem:hitting_lbd}, with probability tending to $1$, $I^N_t\leq \sqrt{I_0m\paragamma}$ for all $t\geq 0$. 

Again, let $M =(M_t)_{t\geq 0}$ be the immigration-death process dominating $R^N$.
Let $M_0 = m(N)$, and we have $k = o(M_0\paragamma)$. For 
\[t_0 = m\paragamma  =\smallo{\left(\frac{M_0\paragamma}{ke}\right)^{M_0}}, \] 
by \Cref{lem:hitting_id_nonabsorb} and the argument similar to the previous case, we have
\[\prob{}{R^N_t\geq 2m(N), \,\forall t\leq t_0\, \labs I^N_t\leq\sqrt{I_0m\paragamma},\,\forall t\geq 0}\to 0.\]

The extinction times $T_{\underline{L}}$ and $T_{\overline{L}}$ have the same asymptotic distribution as specified in \Cref{thm:bdp_extinction_limit} (Case 2 when $I_0=O(1)$ and Case 4 when $I_0\to\infty$). As in the previous case, as $N\to\infty$, 
\[ \prob{}{T_{\overline{L}} < t_0}\to 1.\]

Since 
$I_0\sqrt{I_0m\paragamma} = o(I_0m\paragamma) = O(I_0M_0)$, and $I_0M_0 \leq I_0N^{1/2}\paragamma^{-1/2}+I_0R_0 =  o(N)$,
we have
\begin{align*}
     \lim_{N\to\infty}\left(1-\paralambda\left(1-\frac{\sqrt{I_0m\paragamma}+2m}{N}\right)\right)\underline{L}_0 & = \lim_{N\to\infty}(1-\paralambda)I_0+ \lim_{N\to\infty}\paralambda\frac{I_0(\sqrt{I_0m\paragamma}+2m)}{N}\\
     & =\lim_{N\to\infty}(1-\paralambda)\overline{L}_0 = a.
\end{align*}

\textbf{Case 1.3: }$I_0(1-\paralambda)\to \infty$, $I_0 = \smallo{\frac{N(1-\paralambda)\paragamma}{\log I_0(1-\paralambda)}}$, and $R_0 = \smallo{\frac{N(1-\paralambda)}{\log I_0(1-\paralambda)}}$.

This case is possible only if \eqref{eq:para_sub1} is true.
It covers the scenarios where $\paralambda$ is independent of $N$, and $I_0 \to \infty$.

Let $k(N) := 2I_0$.
Since $ (1-\paralambda)I_0\to\infty$ and $\paralambda(N)< 1$, by the first case in \Cref{lem:hitting_lbd}, with probability tending to $1$, $I^N_t\leq 2I_0$. 

Since $I_0 = \smallo{\frac{N(1-\paralambda)\paragamma}{\log I_0(1-\paralambda)}}$, and $R_0 = \smallo{\frac{N(1-\paralambda)}{\log I_0(1-\paralambda)}}$, we can find $\widetilde{m}(N)$ such that 
\[I_0 \ll \widetilde{m}\paragamma \ll \frac{N(1-\paralambda)\paragamma}{\log I_0(1-\paralambda)}.\]
Let 
\[m(N) = \frac{N(1-\paralambda)}{\log^2 N(1-\paralambda)}\vee R_0\vee \widetilde{m}.\]
We have the properties:
$I_0=o(m\paragamma)$ and 
\[m = \smallo{\frac{N(1-\paralambda)}{\log I_0(1-\paralambda)}}.\] 
By \eqref{eq:para_sub1}, we also have $(1-\paralambda)^{-1}\ll N(1-\paralambda)$.

Define linear birth-death processes $\overline{L}$ and $\underline{L}$ the same way as in Case 1.2.

The extinction time of $T_{\overline{L}}$, according to Case 5, \Cref{thm:bdp_extinction_limit}, is of order $(1-\paralambda)^{-1}\log I_0(1-\paralambda)$. Notice that  
\[(1-\paralambda)^{-1}\log I_0(1-\paralambda)\ll N^2(1-\paralambda)^{2}.\]
Let $t_0 = N^2(1-\paralambda)^2$, then similarly, we have 
\[ \prob{}{T_{\overline{L}} < t_0}\to 1.\]
Since
\[\log t_0  = 2\log N(1-\paralambda) \ll N^{1/2}(1-\paralambda)^{1/2}\log\frac{m\paragamma}{2eI_0}\ll m\log \frac{m\paragamma}{2eI_0},\]
it follows from \Cref{lem:hitting_id_nonabsorb} that,
\[\prob{}{R^N_t\geq 2m(N), \,\forall t\leq t_0\, |I^N_t\leq 2I_0,\,\forall t\geq 0}\to 0.\]

Since $I_0=\smallo{m(N)}$, we have
\begin{align*}
     &\lim_{N\to\infty}\left(1-\paralambda\left(1-\frac{2(I_0+m)}{N}\right)\right)\underline{L}_0 -\lim_{N\to\infty}(1-\paralambda)I_0 =  \lim_{N\to\infty}\paralambda\frac{2I_0(I_0+m)}{N}\\
     & =\smallo{\frac{I_0(1-\paralambda)}{\log I_0(1-\paralambda)}},
\end{align*}
and the rest follows from the order-preserving coupling as introduced in \Cref{appd:coupling_for_bdp}. 

\subsection{Proof of Cases 2.1-2.2}
For Cases 2.1 and 2.2, we require $\paragamma$ to be sufficiently small, so that $R^N$ does not move far away from $R_0(N) \sim r_0 N$, $r_0\in(0,1)$ before extinction. According to \Cref{thm:bdp_extinction_limit}, when $I_0 = O(1)$, we expect the extinction time to be of order $O(1)$; whereas when $I_0\to\infty$, we expect the extinction time to have the asymptotic expansion $\log I_0 + O(1)$.\\

Firstly, we estimate the probability that $R^N$ will remain close to $r_0N$ for a duration of order $\log N$. The approach we use below is a variation of the ODE approximation of Markov chains. A comprehensive introduction to this can be found in \citet{darling2008}. We state the proposition in \Cref{appd:ode_limit} related to our proof.
\begin{lem}\label{lem:R_near_r0N}
Let 
\[X^{N,1}_t := I^N_t/N,\quad X^{N,2}_t := R^N_t/N,\]
with initial states $X^{N,1}_0 = I_0(N)/N$ and $X^{N,2}_0 \to r_0>0$. Let $\delta = \delta(N)>0$. For sufficiently large $N$, if $t_1 = t_1(N)$ satisfies $0<t_1 < \delta\paragamma^{-1}$, then for such $N$ we have
\begin{align}\label{eq:near_bound_p2}
    \prob{}{ \sup_{t\leq t_1}\abs{X^{N,2}_t - X^{N,2}_0}>4\delta}\leq 2\exp{-\frac{\delta^2 N}{4(\paragamma+1)t_1}}+\frac{I_0}{(1-\paralambda+\paralambda r_0/2)\delta N}.
\end{align}
\end{lem}
\begin{proof}
We consider $N\in\N$ to be sufficiently large and fixed throughout the proof.

For $\epsilon>0$, let 
\[T_R(\epsilon):=\inf\left\lbrace t\geq 0: \sup_{s\leq t} \abs{X^{N,2}_s-X^{N,2}_0}>\epsilon \right\rbrace.\]
The process $X^N$ has transition rates:
\begin{align}\label{def:RN_q_j}
q^N((x_1,x_2) ,j)= 
\begin{cases}
       & N\paragamma x_2, \quad j = (0,-\frac{1}{N}), \\
       & Nx_1,\quad j = (-\frac{1}{N},\frac{1}{N}),\\
       & N\paralambda(1-x_1-x_2)x_1,  \quad j = (\frac{1}{N},0) .
\end{cases}
\end{align}
It is also easy to see that the state space of $X^N$ is a subset of $[0,1]^2$.\\

By the argument introduced in \Cref{appd:ode_limit}, we can write 
\begin{align}\label{eq:near_decompose_R}
X^{N,2}_t & = X^{N,2}_0+\int_0^t \sum_j j_2 q^N\left((X^{N,1}_s, X^{N,2}_s) ,j\right)ds + M^N_t \nonumber\\
& = X^{N,2}_0+\int_0^t \left(-\paragamma X^{N,2}_s+X^{N,1}_s\right)ds+ M^N_t,
\end{align}
where $M^N$ is a zero-mean martingale.
We also have for any $x_1,x_2\in[0,1]$,
\[\sum_{j\in J^N} j_2^2q^N((x_1, x_2), j) = N^{-1}\paragamma x_2 + N^{-1}x_1 < (\paragamma+1)N^{-1},\]
where $j_i$ denotes the $i$-th component of $j$.

For any given $N$ and $\epsilon>0$, let
\[T_M(\epsilon):=\inf\left\lbrace t: \sup_{s\leq t}\abs{ M^N_s} >\epsilon \right\rbrace.\]
By \Cref{lem:martingale_bound}, we have for any $t_1=t_1(N)$, 
\begin{align}\label{eq:near_T_W}
 \prob{}{T_M(\epsilon) \leq t_1 } \leq 2\exp{-\frac{\epsilon^2 N}{4(\paragamma+1)t_1}}.   
\end{align}
Taking the supremum and applying Gronwall's inequality to \eqref{eq:near_decompose_R}, we have
\begin{align}
& \sup_{s\leq t}X^{N,2}_s\leq X^{N,2}_0+\int_0^t \paragamma \sup_{u\leq s}X^{N,2}_u ds + \int_0^t X^{N,1}_s ds + \sup_{s\leq t}\abs{M^N_s},\nonumber\\
& \sup_{s\leq t}X^{N,2}_s\leq \left( X^{N,2}_0+\sup_{s\leq t}\abs{M^N_s}+ \int_0^t X^{N,1}_s  ds\right) e^{\paragamma t},\nonumber\\
& \sup_{s\leq t}\left(X^{N,2}_s-X^{N,2}_0\right) \leq X^{N,2}_0(e^{\paragamma t}-1)+\left(\sup_{s\leq t}\abs{M^N_s}+ \int_0^t X^{N,1}_s ds\right) e^{\paragamma t}. \label{eq:near_sup_R}
\end{align}

On the other hand, from \eqref{eq:near_decompose_R} we have,
for all $t\geq 0$, 
$X^{N,2}_t\geq X^{N,2}_0-\paragamma t -\sup_{s\leq t}\abs{M^N_s}$. 
It follows that for all $t>0$, 
\begin{align}\label{eq:near_inf_R}
\inf_{s\leq t}\left(X^{N,2}_s-X^{N,2}_0\right) \geq -\paragamma t- \sup_{s\leq t}\abs{M^N_s}.
\end{align}
Combining \eqref{eq:near_sup_R} and \eqref{eq:near_inf_R}, we have
\[
\sup_{s\leq t}\abs{ X^{N,2}_s-X^{N,2}_0} \leq (e^{\paragamma t}-1)+\left(\sup_{s\leq t}\abs{M^N_s}+ \int_0^t X^{N,1}_sds\right) e^{\paragamma t}. 
\]

Define $T_{int}(\epsilon):=\inf\left\lbrace t:\int_0^t X^{N,1}_sds>\epsilon \right\rbrace$.

For $t_1(N)$ and $\delta = \delta(N)\to 0$ satisfying $\paragamma(N)t_1(N)<\delta(N)$ for sufficiently large $N$, on the event 
\[\left\lbrace t< T_M(\delta)\wedge T_{int}(\delta)\wedge t_1\right\rbrace ,\] 
we have
\[
\sup_{s\leq t}\abs{ X^{N,2}_s-X^{N,2}_0}  = \left(\paragamma t_1+ O(\paragamma^2 t_1^2)\right)+2\delta(1+\paragamma t_1+O(\paragamma^2 t_1^2)) < 4\delta.
\]
In other words, $\prob{}{ T_R(4\delta)> t\labs t< T_M(\delta)\wedge T_{int}(\delta)\wedge t_1}=1$. 

It follows that 
\begin{align}\label{eq:near_T_p2}
    \prob{}{ T_M(\delta)\wedge T_{int}(\delta)\wedge t_1\leq T_{R}(4\delta)}=1.
\end{align}
Let $T_{sum}(x):=\inf\left\lbrace t:X^{N,1}_t+X^{N,2}_t\leq x \right\rbrace$.

It follows from \eqref{eq:near_T_p2} that for sufficiently large $N$, on the event $\left\lbrace t<T_M(\delta)\wedge T_{int}(\delta)\wedge t_1\right\rbrace $,
\[X^{N,1}_t+X^{N,2}_t\geq X^{N,2}_0-4\delta > \frac{3r_{0}}{4},\]
which suggests that
\begin{align}\label{eq:near_T_X1}
\prob{}{T_{sum}\left(\frac{r_0}{2}\right) > T_M(\delta)\wedge T_{int}(\delta) \wedge t_1} = 1.
\end{align}
The inequality is strict because $(X^{N,1}+X^{N,2})$ has jump sizes of order $N^{-1}$ and cannot reach $r_0/2$ from above $3r_0/4$ in one jump.

Also from \eqref{eq:near_T_p2},
\[\prob{}{ T_{R}(4\delta)< t_1} \leq \prob{}{ T_M(\delta)< t_1}+ \prob{}{T_{int}(\delta) = T_{int}(\delta)\wedge T_M(\delta) < t_1}.\]
The upper bound of $\prob{}{ T_M(\delta)< t_1}$ is obtained in \eqref{eq:near_T_W}. For the second term on the RHS above, 
conditioning on the event $\left\lbrace
T_{int}(\delta) = T_{int}(\delta)\wedge T_M(\delta) < t_1\right\rbrace$, the equality \eqref{eq:near_T_X1} is equivalent to
\[\prob{}{ T_{sum}\left(\frac{r_0}{2}\right) > T_{int}(\delta)} = 1.\]
Then
\[\prob{}{T_{int}(\delta) = T_{int}(\delta)\wedge T_M(\delta) < t_1} \leq \prob{}{T_{int}(\delta) < T_{sum}\left(\frac{r_0}{2}\right)\wedge t_1}.\]

On the event $\left\lbrace t < T_{sum}\left(\frac{r_0}{2}\right)\wedge t_1\right\rbrace$, the process $I^N$ is dominated by a linear birth-death process $L= (L_t)_{t\geq 0}$ with birth rate $\paralambda\left(1-\frac{r_0}{2}\right)$ and death rate $1$. Therefore, $\int_0^tX^{N,1}_s ds$ is stochastically bounded by \[N^{-1}\int_0^{T_L}L_sds = N^{-1}H(I_0),\] 
where $T_L$ and $H(I_0)$ are defined as \Cref{lem:integral_BDP}. \\

For any $t>0$, the probability $\prob{}{ T_{int}(\delta) \leq t}$ is then bounded by the probability  $\prob{}{H(I_0)>N\delta}$. By \eqref{eq:integralBDP_bound}, 
\begin{align*}
    & \prob{}{T_{int}(\delta) = T_{int}(\delta)\wedge T_M(\delta) < t_1} \leq \prob{}{T_{int}(\delta) < t_1}\leq \prob{}{H(I_0)>N\delta}\\
    \leq & \frac{I_0}{(1-\paralambda+\paralambda r_0/2)N\delta},
\end{align*}
where the last inequality follows from \Cref{lem:integral_BDP}.

Together with \eqref{eq:near_T_W}, we have
\[
\prob{}{ T_{R}(4\delta)< t_1} \leq 2\exp{ -\frac{\delta^2 N}{4(\paragamma+1)t_1}}+\frac{I_0}{(1-\paralambda+\paralambda r_0/2)N\delta}.
\]
\end{proof}
Now we are ready to discuss different cases under the second scenario, depending on the size of $I_0$.

\textbf{Case 2.1: } $I_0 = O(1)$ and $\paragamma = o(1)$.\\

Let $\delta(N) = \paragamma^{1-\epsilon}\vee N^{-1/3}$ for a small $\epsilon>0$. Then  for any $t_1 = O(1)$ we have $t_1<\paragamma^{-\epsilon}\leq \delta\paragamma^{-1}$. By \Cref{lem:R_near_r0N},
\[\prob{}{\sup_{s\leq t_1}\abs{R^N_s-R_0}>4\delta N}\to 0.\]

Define $\delta_0(N) := \abs{R_0/N-r_0}$.

For each $N\in\N$, define two linear birth-death processes $
\underline{L}$ and $\overline{L}$ such that $\overline{L}$ has birth rate $\paralambda(1-r_0+4\delta+\delta_0)$ and death rate $1$, and $\underline{L}$ has birth rate $\paralambda(1-r_0-4\delta-\delta_0)$  and death rate 1. Let $\underline{L}_0=\overline{L}_0 = I_0$. 

Denote \[T_R(4\delta):=\inf\left\lbrace t\geq 0: \sup_{s\leq t} \abs{R^N_s-R_{0}}>4\delta N \right\rbrace.\]
Then for $t< T_R(4\delta)$,
\[I_0\paralambda(1-r_0-4\delta-\delta_0)\leq I_0 \paralambda\left(1-\frac{I^N_t+R^N_t}{N}\right)\leq I_0\paralambda(1-r_0+4\delta+\delta_0),\]
and all three terms tend to $\lim_{N\to\infty}I_0\paralambda(1-r_0)$.\\

With $a =\lim_{N\to\infty}\left(1-\paralambda+\paralambda r_0\right)\in (0,\infty)$, the conclusion then follows from Case 2 of \Cref{thm:bdp_extinction_limit}:
\begin{align*}
    \prob{}{T^N_o\leq  w} &  =\prob{}{T^N_o\leq w,\, T_R(4\delta)\leq w}+\prob{}{T^N_o\leq w,\, T_R(4\delta)> w}\\
    & \to\left(1+\frac{a}{e^{aw}-1}\right)^{-I_0},
\end{align*}
where by \eqref{eq:near_bound_p2}, we have for all $w>0$,
\[\prob{}{T^N_o\leq w,\, T_R(4\delta)\leq w}\leq \prob{}{T_R(4\delta)\leq w}\to 0.\]
\textbf{Case 2.2} $I_0\to\infty$ and there exists $\epsilon_1,\epsilon_2>0$ such that $I_0 = \smallo{N^{1-\epsilon_1}}$ and $\paragamma = \smallo{N^{-\epsilon_2}}$.\\

Let $\delta = N^{-b\epsilon_2}$, where positive constant $b$ is chosen to satisfy $b<\frac{1}{2\epsilon_2}\wedge \frac{\epsilon_1}{\epsilon_2}\wedge 1$.
Then  we have that any $t_1 = O(\log N)$ satisfies $t_1<\delta\paragamma^{-1} = N^{(1-b)\epsilon_2}$.
 
 Since $\delta^2N = N^{1-2b\epsilon_2}$ and $\frac{I_0}{\delta N} \ll N^{b\epsilon_2-\epsilon_1}$ are both of negative polynomial orders of $N$, by \Cref{lem:R_near_r0N},
\[\prob{}{\sup_{s\leq t_1}\abs{R^N_s-R_0}>4\delta N}\to 0.\]
The constructions of $\underline{L}$ and $\overline{L}$, and the stopping time $T_R(4\delta)$ remain the same as Case 2.1. We still have that, with probability tending to $1$, for $t\leq T_R(4\delta)$,
\[\paralambda(1-r_0-4\delta-\delta_0)\leq \paralambda\left(1-\frac{I^N_t+R^N_t}{N}\right)\leq \paralambda(1-r_0+4\delta+\delta_0).\]

Recall that $a=\lim_{N\to\infty}\left(1-\paralambda+\paralambda r_0\right)\in (0,\infty)$. 
Following from Case 5 of \Cref{thm:bdp_extinction_limit}, the extinction times of $\underline{L}$ and $\overline{L}$ tend to the same limit if $I_0(4\delta+\delta_0) = \smallo{\frac{aI_0}{\log(aI_0)}}$, which is equivalent to
$N^{-b\epsilon_2}\log N \to 0$.

The conclusion then follows from Case 5 of \Cref{thm:bdp_extinction_limit}:
\begin{align*}
     & \prob{}{aT^N_o -\log(aI_0)\leq w} = \prob{}{T^N_o\leq a^{-1}\log(aI_0)+ a^{-1}w < T_R(4\delta)}\\
    & +\prob{}{T^N_o\leq a^{-1}\log(aI_0)+ a^{-1}w ,\,T_R(4\delta)\leq a^{-1}\log(aI_0)+ a^{-1}w }\\
     \to & e^{-e^{-w}},
\end{align*}
where we use the fact that for any constants $c_1,c_2>0$,
\[\prob{}{T_R(4\delta)\leq c_1\log N+ c_2 }\to 0.\]

\section{Numerical analysis}
We compare our analytic results in \Cref{thm:small_ini_main_result} to the data obtained through Monte Carlo simulations. Two simulation algorithms, SSA and modified $\tau$-leaping, are used to simulate the extinction time. The detail of both algorithms and the motivation of developing the modified $\tau$-leaping method as an approximation to the SSA method can be found in \citet{cao2005avoiding}. Roughly speaking, the classic SSA method is time-consuming when simulating stochastic epidemic models with large populations, and even more so when the model is near-critical. In comparison, the modified $\tau$-leaping method is more efficient. 

Our numerical analysis is implemented in MATLAB(R2019b). For each case in \Cref{thm:small_ini_main_result}, we choose one `parameter-initial state' combination and run 700 simulations for a set of $N$ of different orders (up to $N=10^7$ for the SSA method and up to $N = 10^{12}$ for the modified $\tau$-leaping method. The choice of maximum $N$ is made based on both the computation time and the observed speed of convergence). We set the parameters of the modified $\tau$-leaping to be $n_c = 200$ and $\epsilon = 0.02$.

For each case, the results are presented in three sub-figures. We present the first example in large figures so that the reader can examine the details. We can see that all our asymptotic results provide fairly good approximations.
\bigthreefig{Case 1.1  ($\paralambda(N)<1$): $\paralambda = 1-N^{-1/2}, \, \paragamma = N^{-1/6}$, $I_0 = N^{1/4},\, R_0 = N^{1/2}$.\\
    In all figures, the lines are colour-coded by a gradient from dark red to yellow. A line with a colour closer to yellow  corresponds to the result for a larger $N$. In all of our figures, the dashed lines represent simulations by the modified $\tau$-leaping method and the solid lines represent simulations by the SSA method. The time axes are always scaled according to the scaling of the asymptotic distribution. \\
     Figure (a) presents a randomly chosen sample path from the simulation for each $N$ in log scale, i.e., $\log_N(I^N_t)$ (in the thicker lines) and  $\log_N(R^N_t)$ (in the thinner lines) over the scaled time.\\
     Figure (b) presents the histogram of extinction times for different $N$, normalised so that the sum of the bar areas is less than or equal to 1 (i.e., Figure (b) is the simulated probability density function of the extinction time). The blue line represents the first-order derivative of the asymptotic distribution function.\\
     Figure (c) presents the histogram of extinction times for different $N$, normalised so that the height of the last bar is less than or equal to 1 (i.e., Figure (c) is the simulated cumulative distribution function of the extinction time). The blue line represents the asymptotic distribution we have derived through analysis.
}
{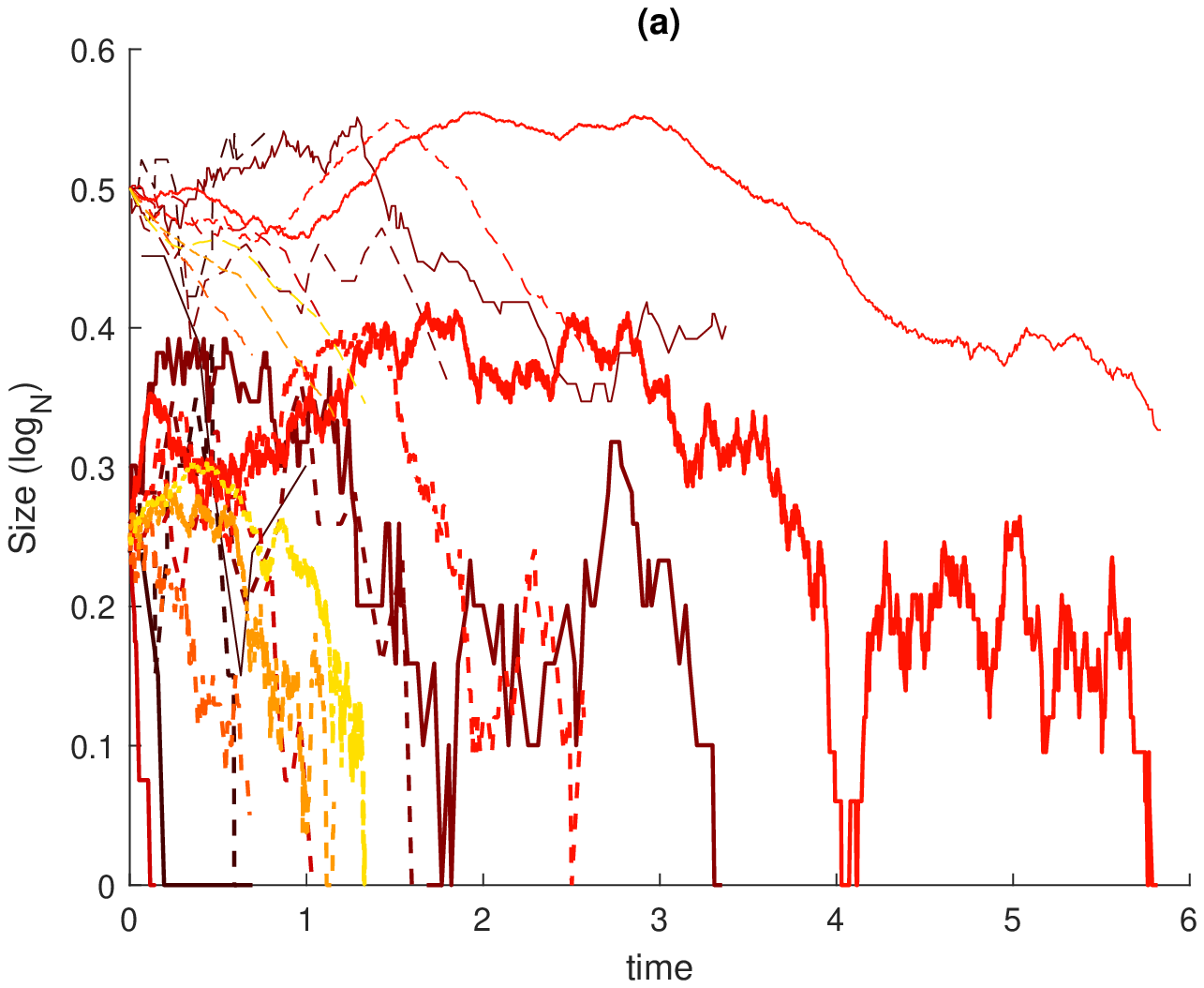}{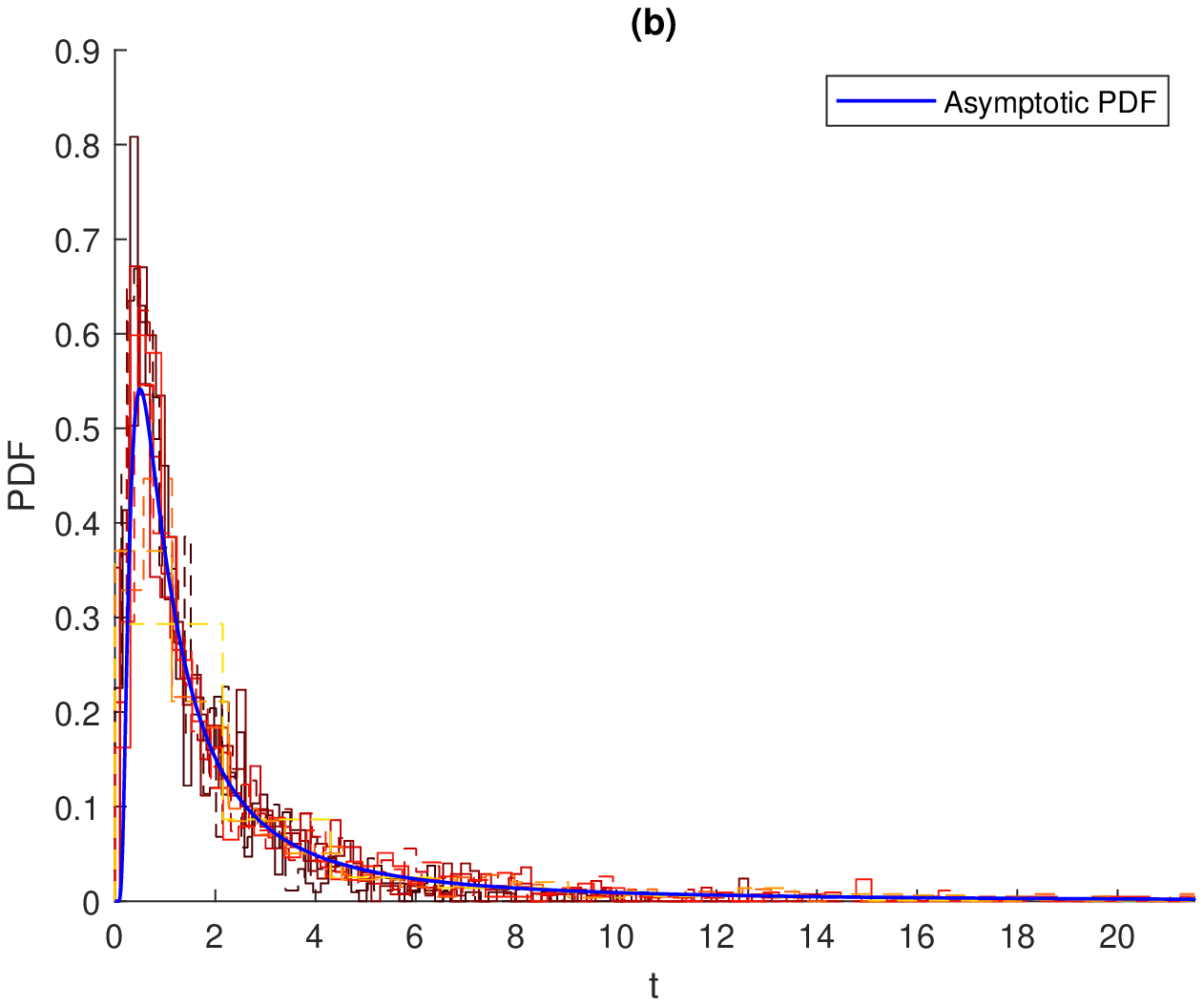}
{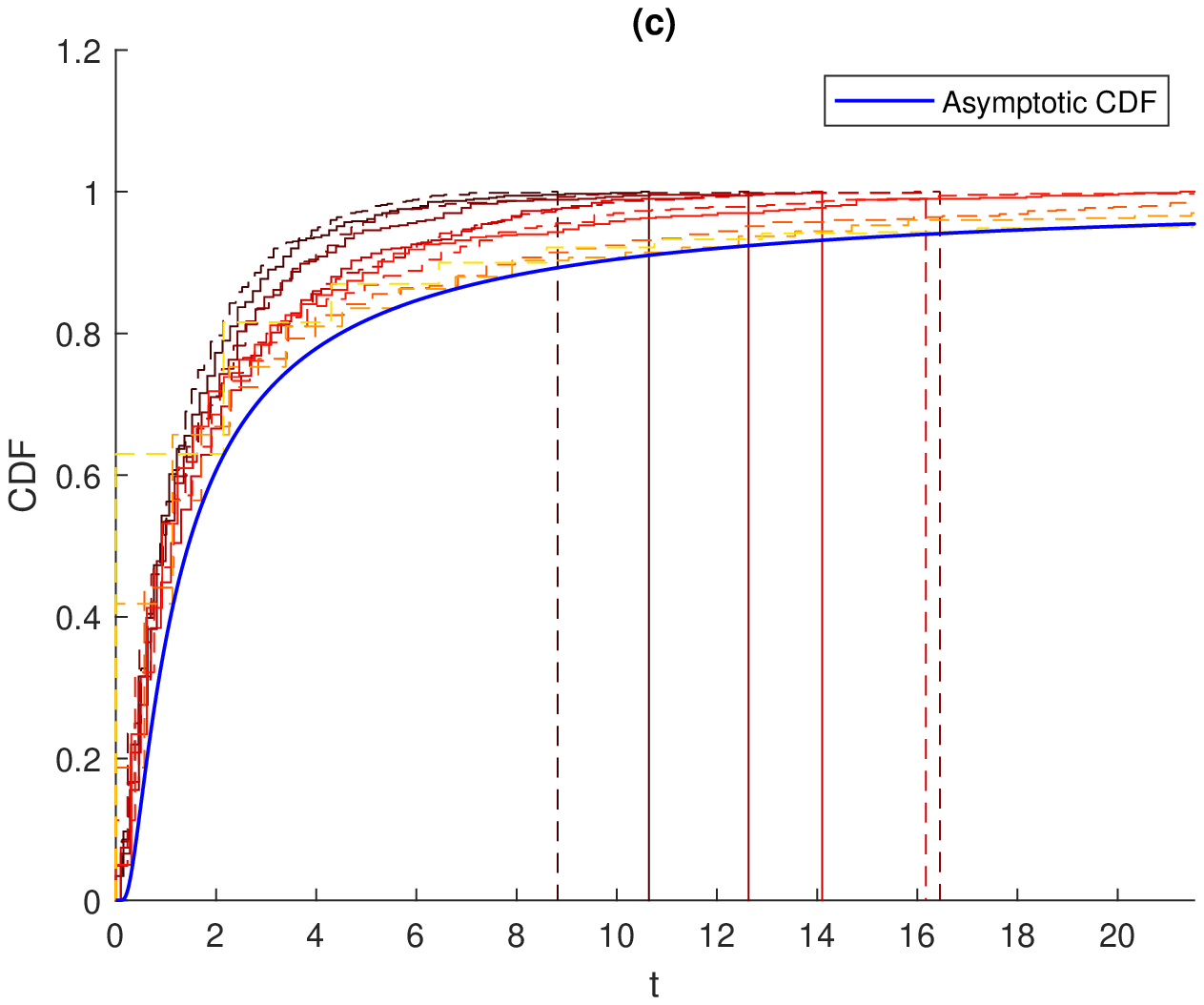}
{fig:case11}
\twobytwofig{
Case 1.1 ($\paralambda(N)>1$): $\paralambda = 1 + N^{-1/2}, \, \paragamma = N^{-1/6}$,  $I_0 = N^{1/4},\, R_0 = N^{1/2}$.
}
{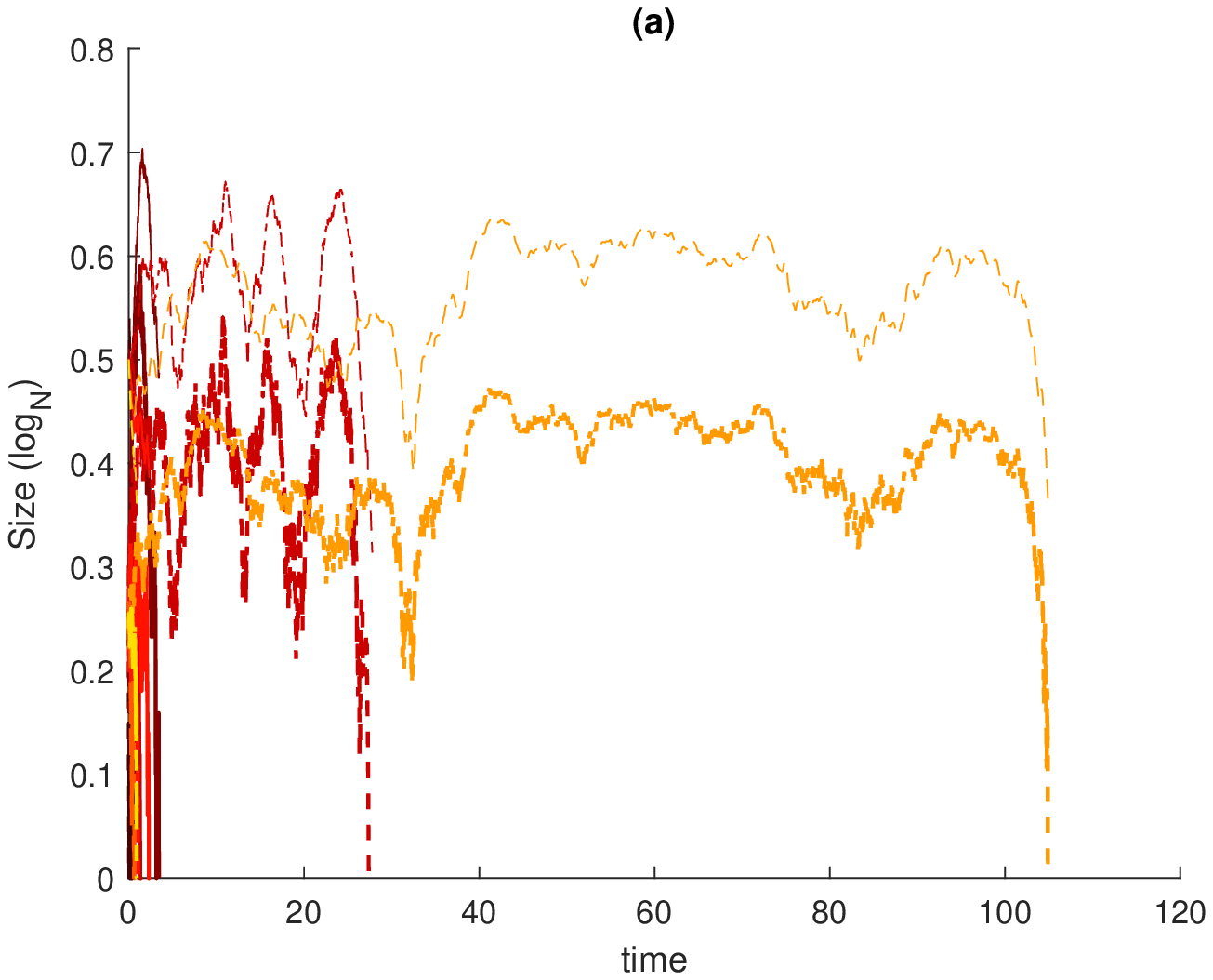}{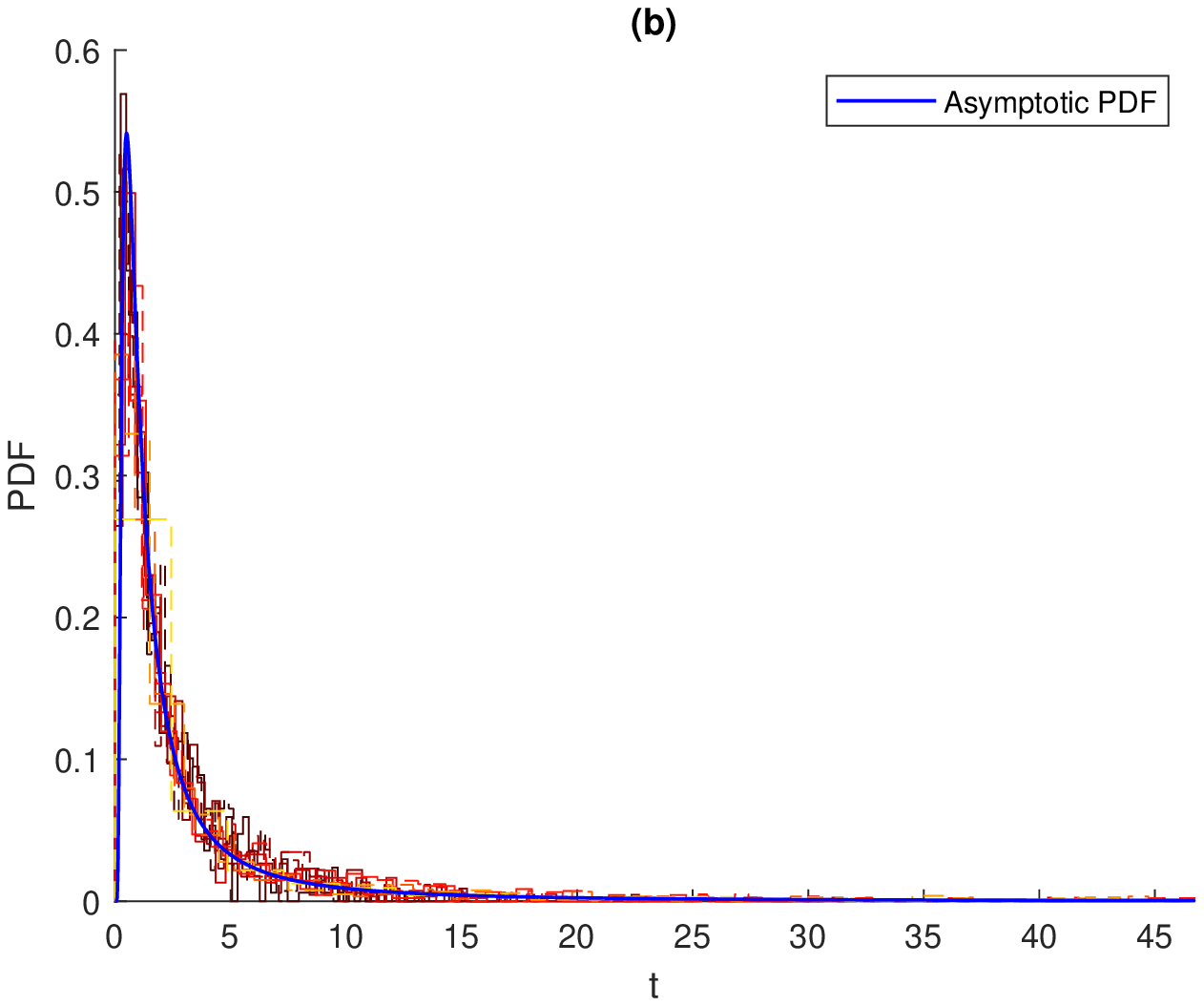}
{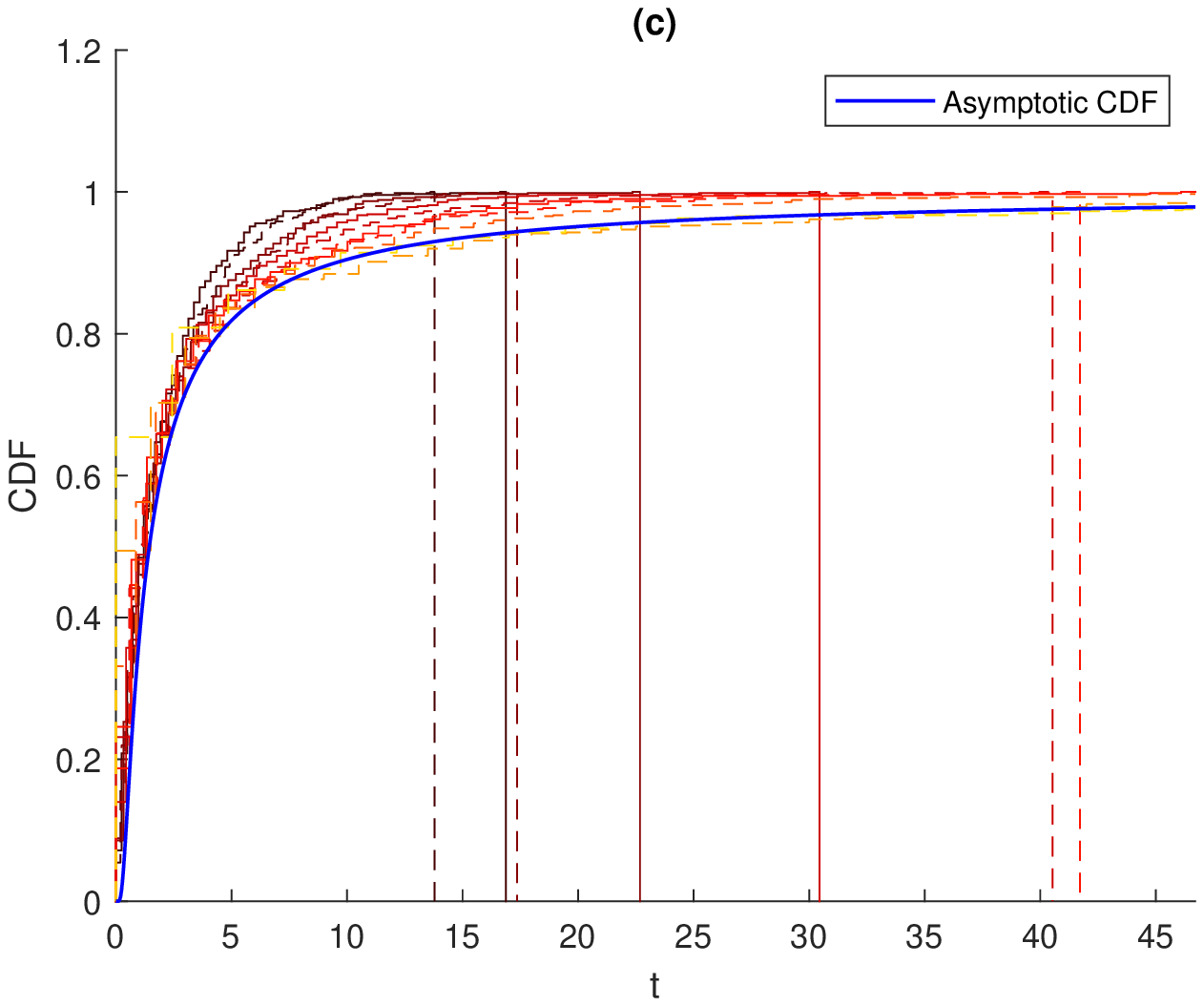}
{fig:case11plus}
\twobytwofig{
Case 1.2: 
    $\paralambda = 1-N^{-1/3}, \, \paragamma = N^{-1/4}$, 
    $I_0 = N^{1/3},\, R_0 = N^{1/3}$.
}
{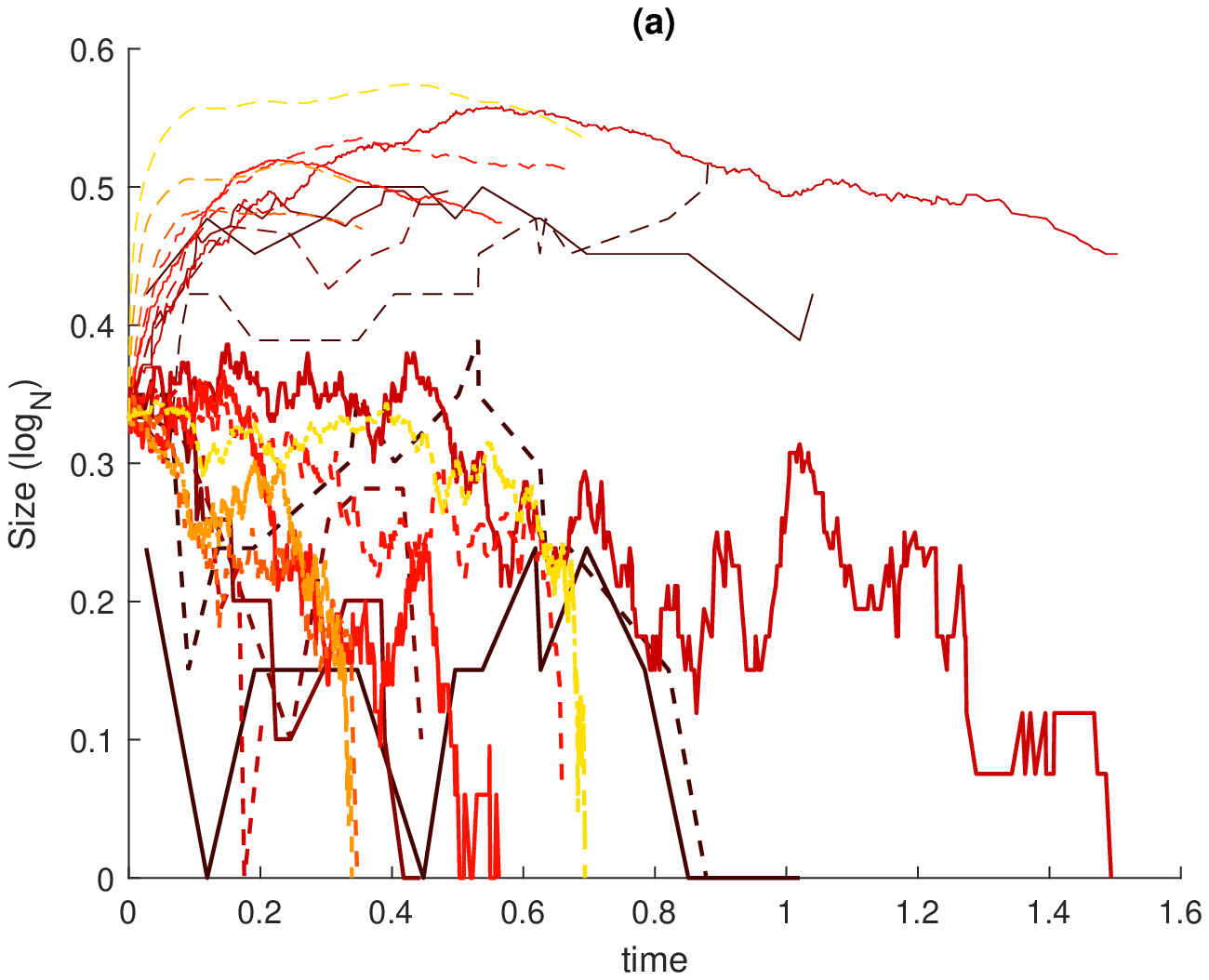}{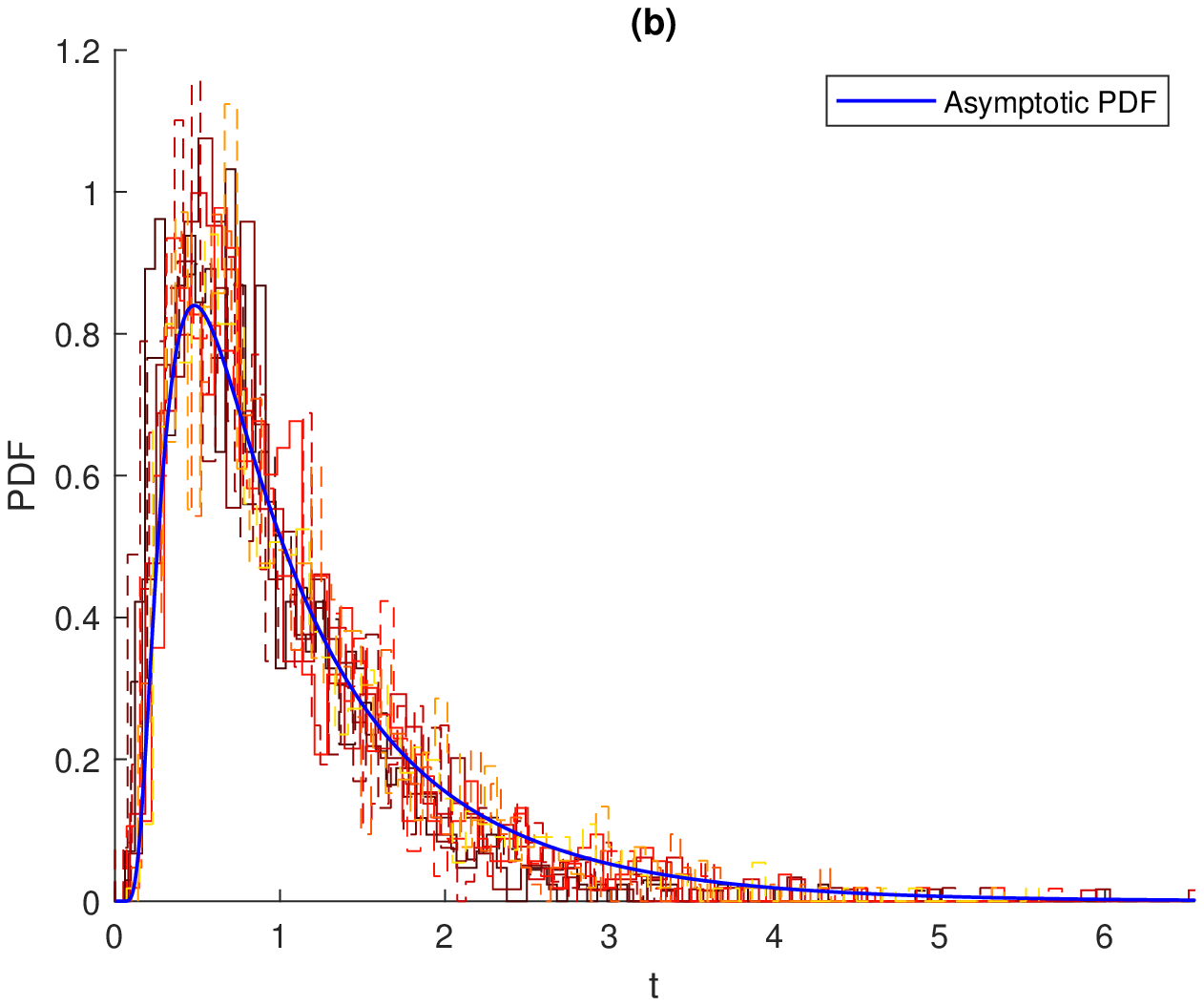}
{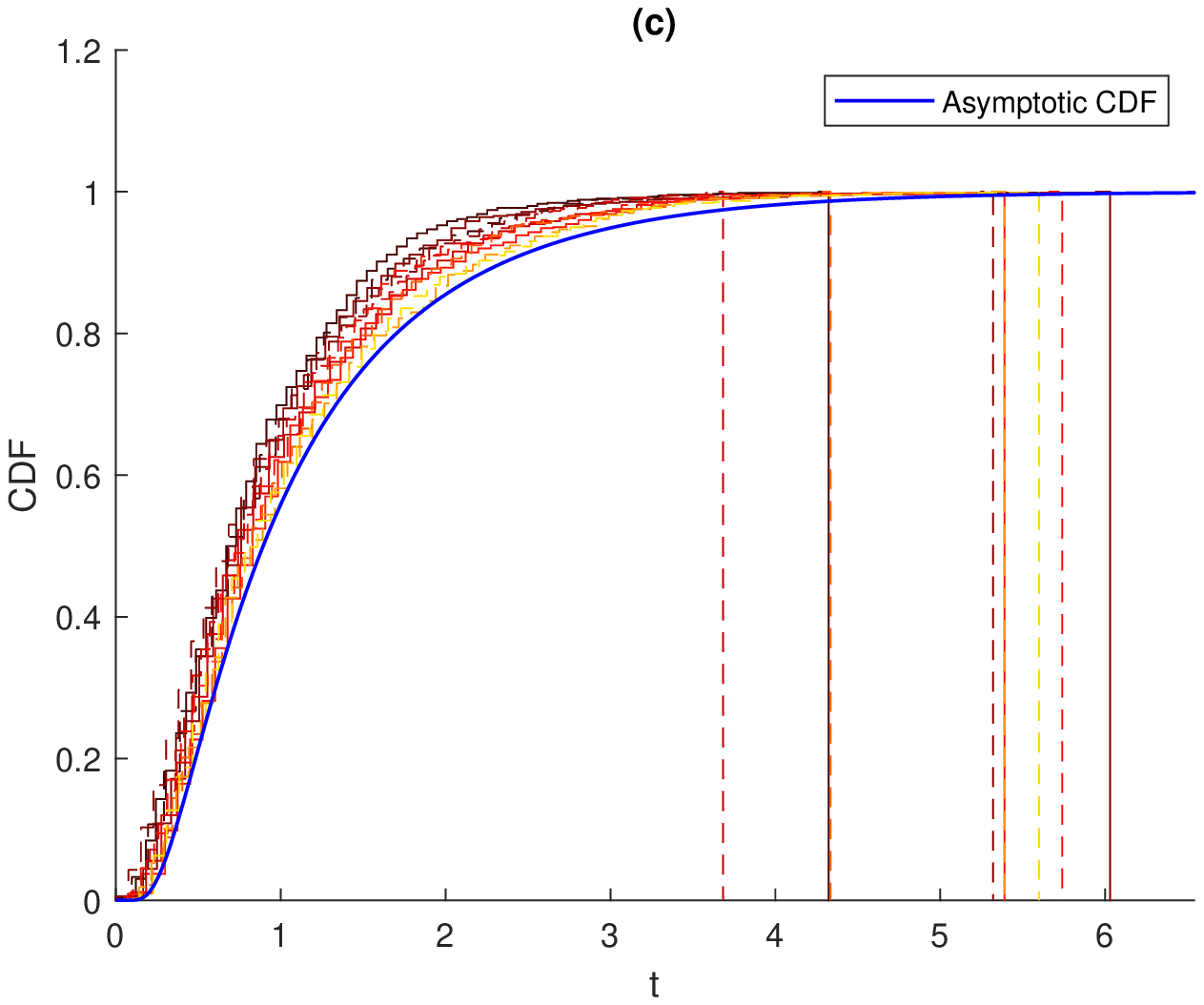}
{fig:case12}
\twobytwofig{
Case 1.3:
$ \paralambda =1- N^{-1/4}, \, \paragamma = N^{-1/4}$, $I_0 = N^{1/3},\, R_0 = N^{1/2}$.
}
{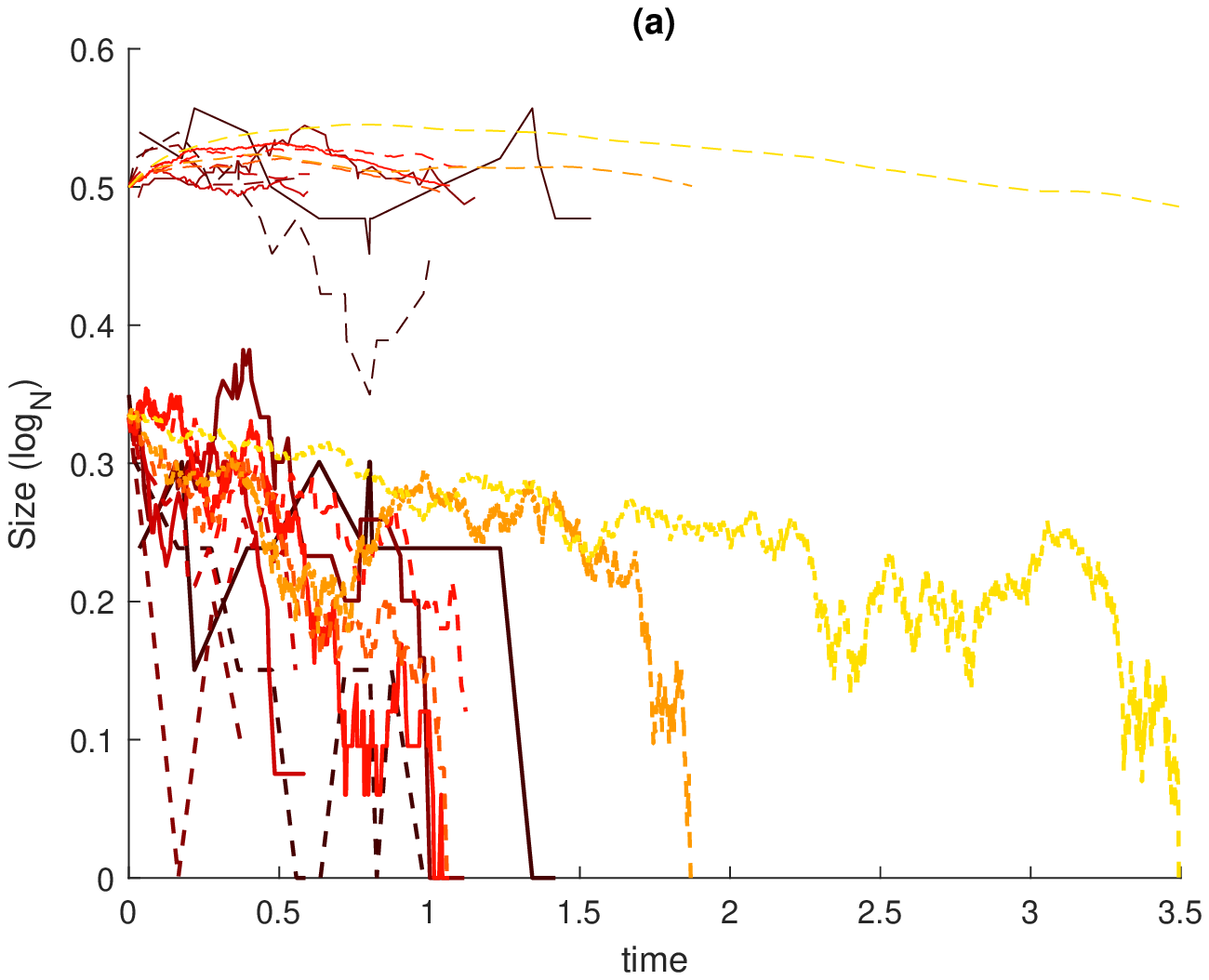}{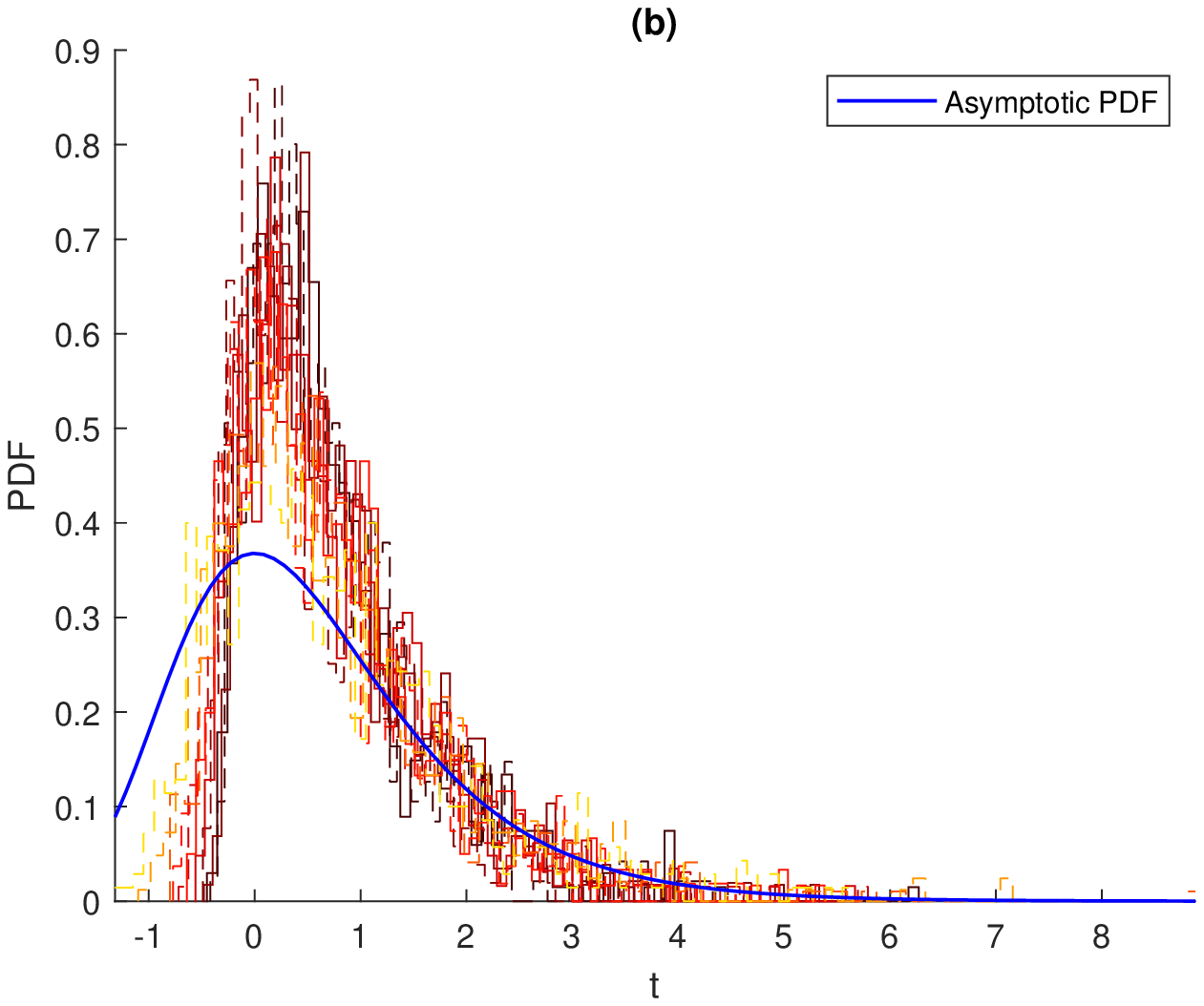}
{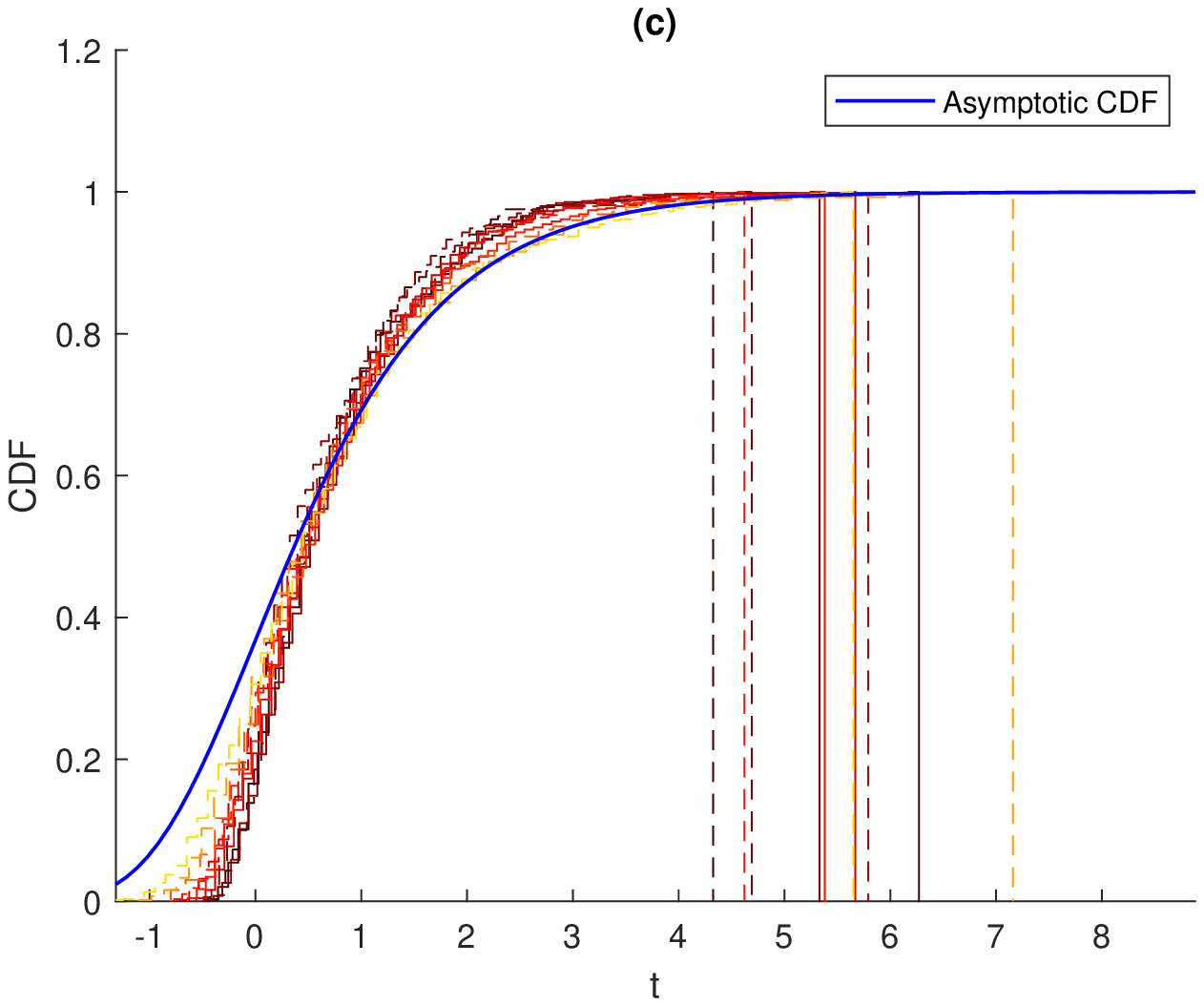}
{fig:case13}

\twobytwofig{
Case 2.1: $\paralambda = 0.8-N^{-1/2}, \, \paragamma = N^{-5/12}$, $I_0 = 30,\, R_0 = 0.3N$.
}
{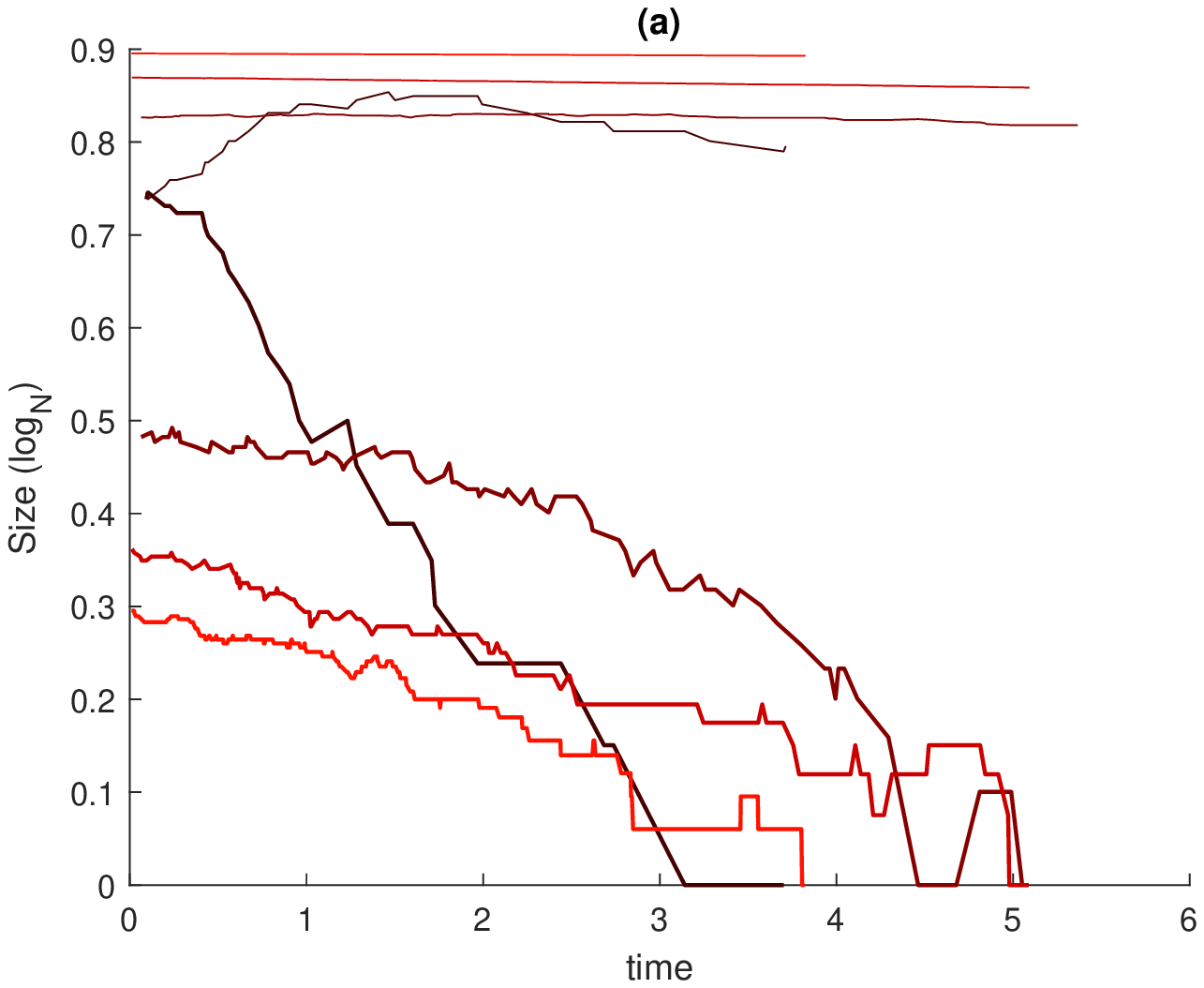}{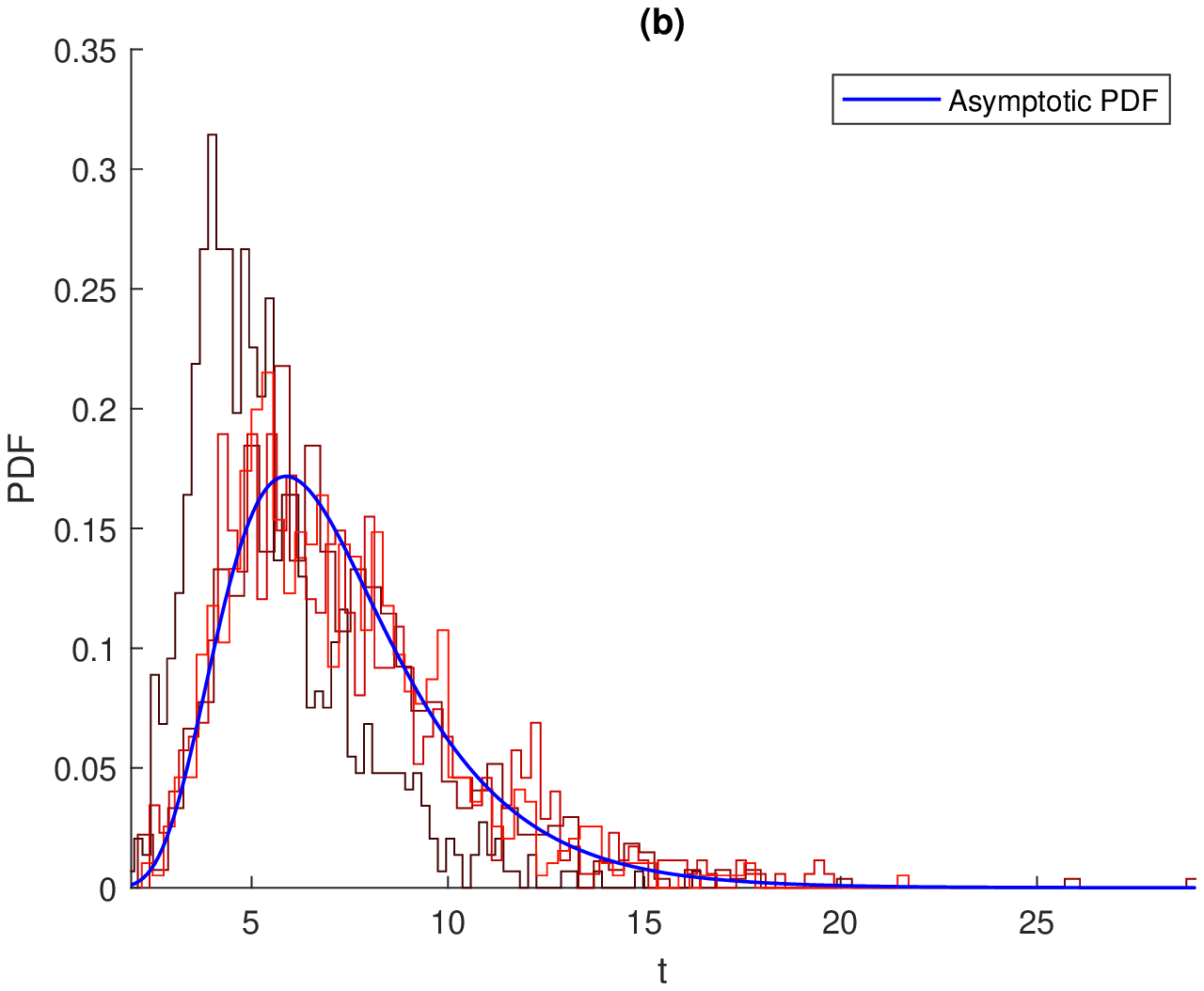}
{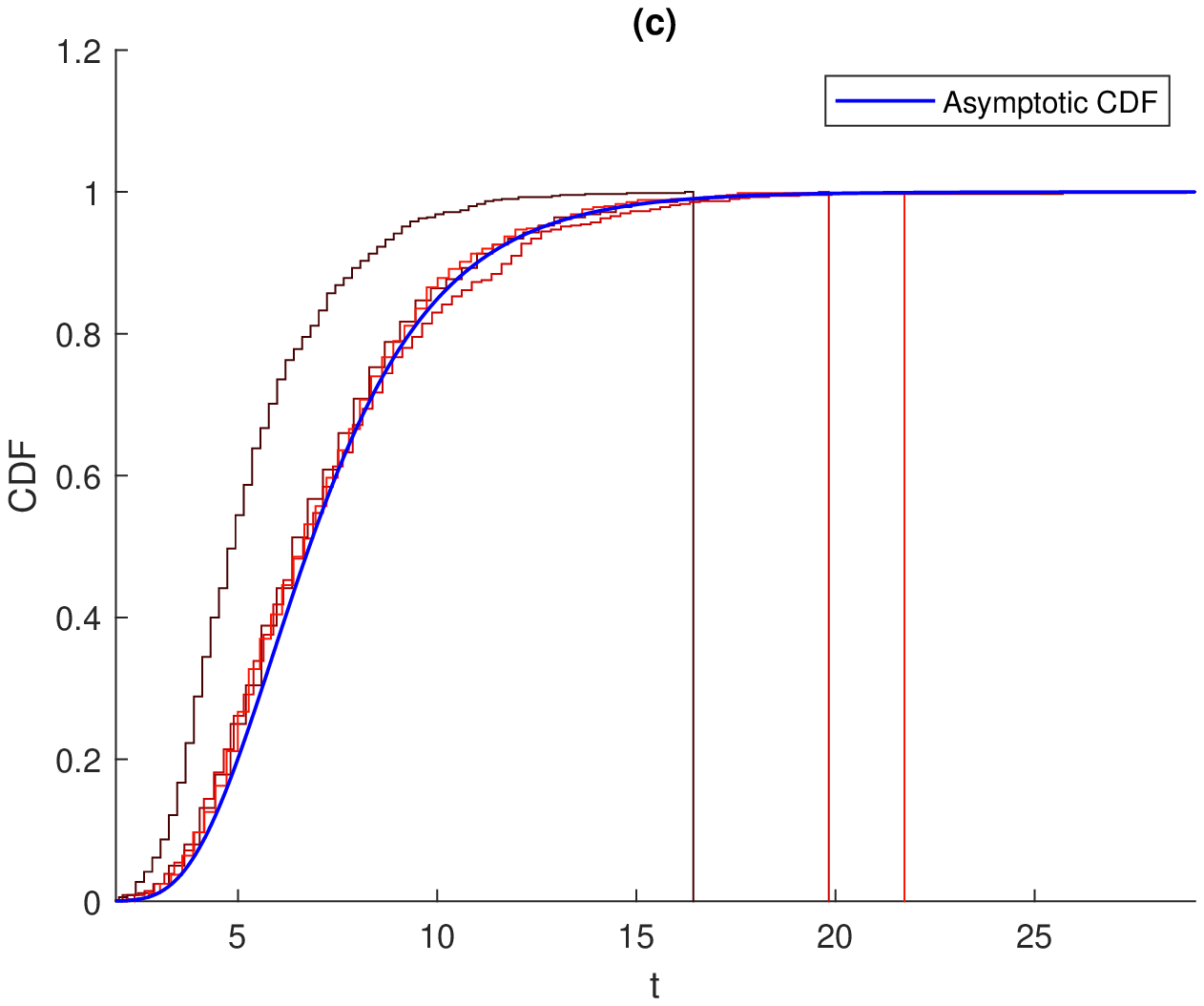}
{fig:case21}

\twobytwofig{
Case 2.2: $\paralambda = 1-N^{-1/4}, \, \paragamma = N^{-1/2}$, $I_0 = N^{1/5},\, R_0 = 0.7N$.
}
{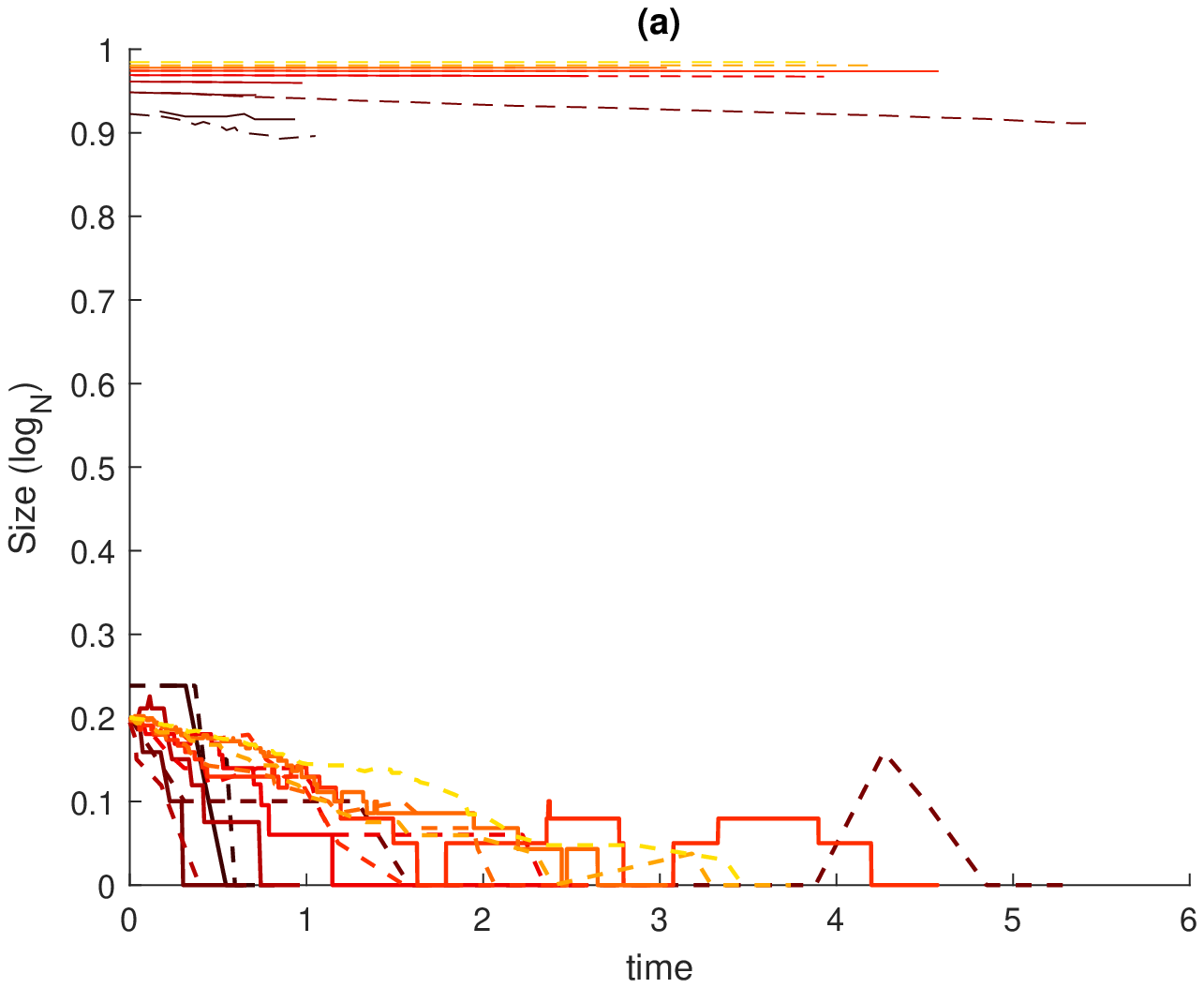}{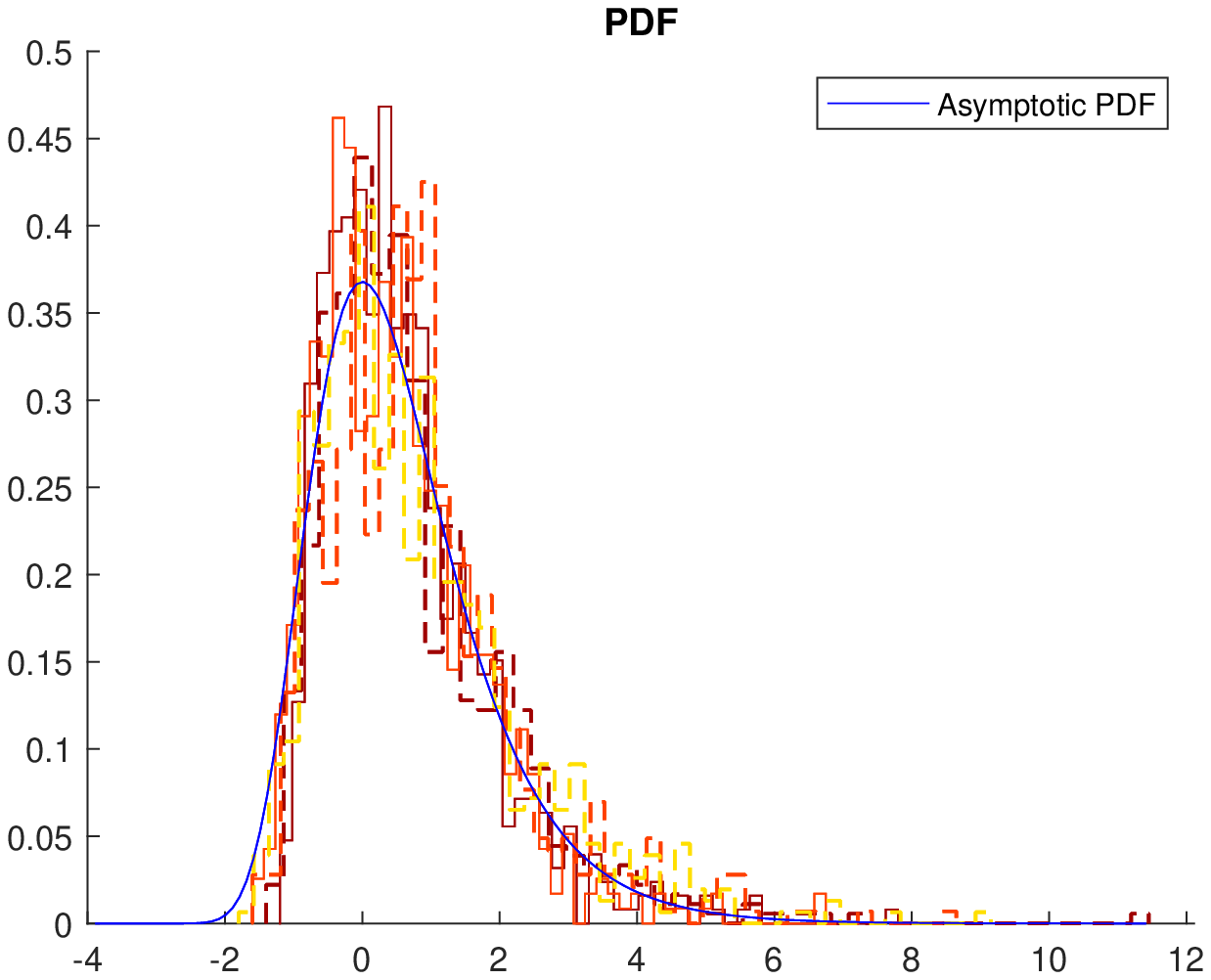}
{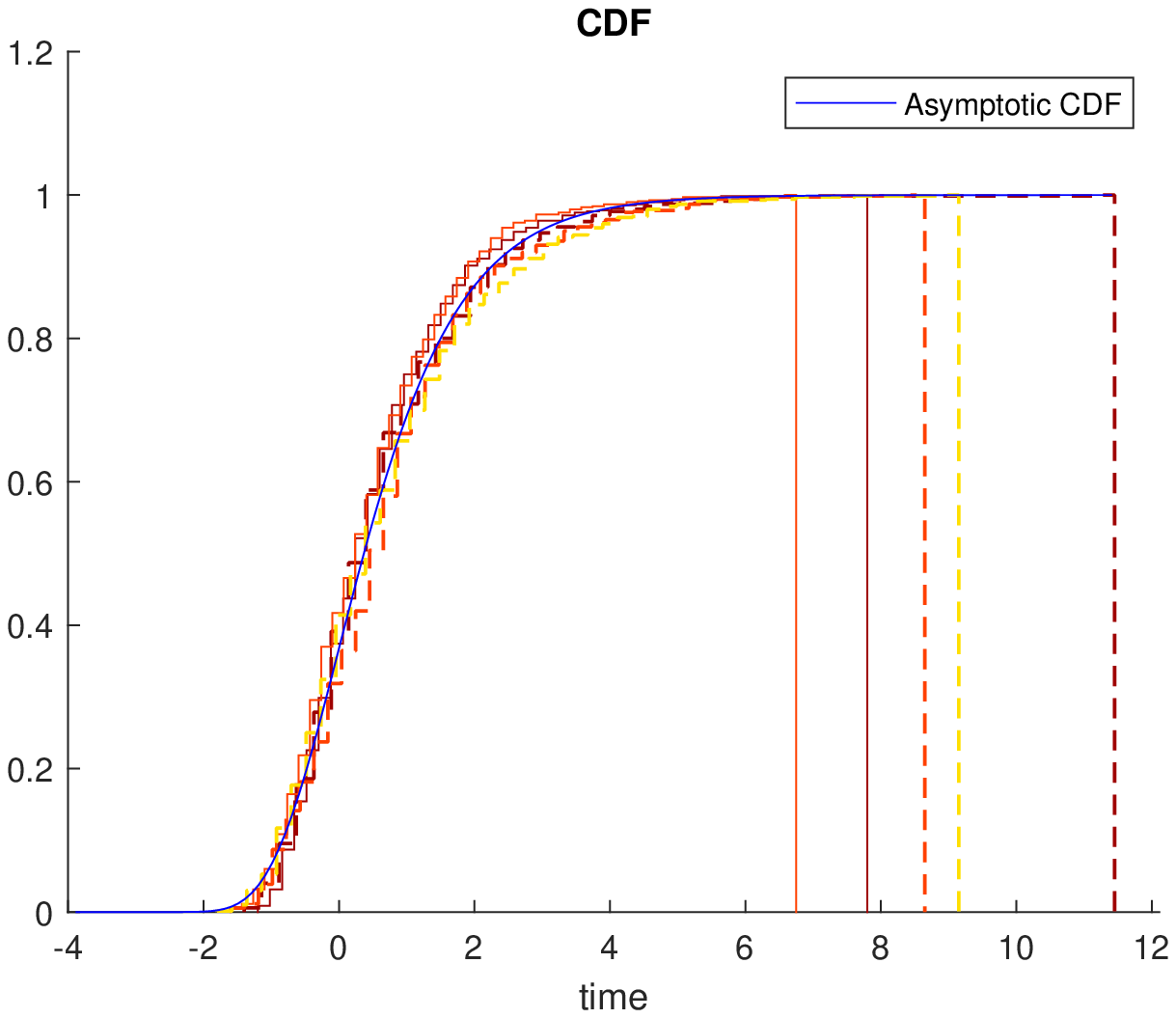}
{fig:case22}

\appendix
\section{Order-preserving coupling for birth-death processes}\label{appd:coupling_for_bdp}
Consider two birth-death processes $Z^{1},Z^{2}$, both defined on a state space $\N$. For $i=1,2$, $Z^{i}$ has transition rates
\[
\begin{cases}
&z \to z+1 \textnormal{, at rate } b_i(z),\\
&z \to z-1 \textnormal{, at rate } d_i(z),\\
\end{cases}
\]
where $b_1(z)\geq b_2(z)$, $d_1(z)\leq d_2(z)$ for all $z\in\N$, and their initial states satisfy $Z^{1}_0\geq Z^{2}_0$. 
The order-preserving coupling $(\hat{Z}^{1},\hat{Z}^{2})$ is defined on $\N\times\N$, satisfying:

At state $(z,z)$, $(\hat{Z}^{1},\hat{Z}^{2})$ has transition rates
\[
\begin{cases}
    &(z,z) \to (z+1,z) \textnormal{, at rate } b_1(z)-b_2(z),\\
    &(z,z) \to (z+1,z+1) \textnormal{, at rate } b_2(z),\\
    &(z,z) \to (z-1,z-1) \textnormal{, at rate } d_1(z),\\
    &(z,z) \to (z,z-1) \textnormal{, at rate } d_2(z)-d_1(z),\\
\end{cases}
\]
and at state $(z_1, z_2)$, $z_1\neq z_2$, $\hat{Z}^{1}$ and $\hat{Z}^{2}$ jump independently. Since $\hat{Z}^{1}$ and $\hat{Z}^{2}$ will a.s. not jump at the same time, their paths will a.s. not cross each other when $|z_1-z_2|=1$.

The intuition as to why this coupling is order-preserving is that, $\hat{Z}^{1}$, with higher birth rates and lower death rates, will stay above $\hat{Z}^{2}$ until they meet, in which case, they will jump together until either $\hat{Z}^{1}$ moves upward, or $\hat{Z}^{2}$ moves downward.

\section{ODE approximation of Markov Chains}\label{appd:ode_limit}

 Consider a sequence of Markov chains indexed by $N$, valued in finite state space $S^N\subset \R^d$, and is denoted as $\{(X^N_t)_{t\geq 0}\}_{N\in\N}$. For each $N$, $X^N$ is uniquely defined by its initial state $X^N_0 = x^N_0\to x_0$ as $N\to \infty$, and transition rates $q^N(x,j),\,j\in J^N$, where $J^N$ is the set of possible jumps in column vectors. We assume the number of elements in $J^N$ is finite and independent of $N$, and 
\begin{align}
    \bar{j}^N:=\max\left\{|j|:j\in J^N\right\} \to 0,\quad N\to\infty.\label{assum:jump_size}
\end{align}
For each $N$, denote a jump at time $t$ as $\Delta X^N_t :=X^N_t-X^N_{t-}$, then we can define random measures on $[0,\infty)\times J$ as 
\[\mu:=\sum_{t:\Delta X^N_t\neq 0} \delta_{(t,\Delta X^N_t)}\text{ and }\upsilon(dt,j) := q(X^N_{t-},j)dt.\] 
Thus we can write 
\[X_t = X_0+ \int_0^t\sum_{j\in J} j \mu(ds,j).\]

Let $(H(t,j))_{t\geq 0}$ be a left-continuous adapted process for each $j\in J$. It is known that if $H$ satisfies for all $t\geq 0$,
\begin{align*}
    \expect{}{\int_0^t\sum_{j\in J} \abs{H(s,j)}\upsilon(ds,j)}<\infty,
\end{align*}{}
then it is known that 
\begin{align*}
   \int_0^t\sum_{j\in J} H(s,j)(\mu-\upsilon)(ds,j)
\end{align*}{}
is a well-defined martingale (see e.g. Theorem 8.4, \cite{darling2008}).

It follows that we can decompose $X$ as 
\begin{align}\label{eq:decompose_X}
    X^N_t = x_0+\int_0^t\sum_{j\in J} jq(X^N_{s-},j)ds + M^N_t,
\end{align}{}
where the martingale part $M^N_t:=\int_0^t\sum_{j\in J} j(\mu-\upsilon)(ds,j)$ is called the \textit{compensated martingale}, and $\int_0^t\sum_{j\in J} jq(X^N_{s-},j)ds$ is called the \textit{compensator}.

\begin{prop}\label{lem:martingale_bound}
    Consider the sequence of Markov chains $\{X^N\}_{N\in\N}$ defined as above. Denote the i-th component of vector $j$ as $j_i$.
    For each given $N$, $t_0>0$ and $\bar{a} = \bar{a}(N)>0$, let 
    \begin{align}\label{def:stoptime_quadratic_drift}
        \tau^{t_0}(i,\bar{a}):=\inf\bigg\{t:\int_0^{t\wedge t_0}\sum_{j\in J^N} q^{N}(X^N_{s-},j)j_i^2ds > \bar{a}(N)\bigg\}, \quad i=1,2,\cdots, d,
    \end{align}
    then for any $\delta = \delta(N) < \bar{a}/\max_{j\in J^N}\abs{j_i}$,
    \begin{align*}
       \prob{}{\sup_{s\leq t_0\wedge  \tau^{t_0}(i,\bar{a})}\abs{\product{e_i,M^N_s}}\geq \delta}\leq 2\exp{-\frac{\delta^2}{4\bar{a}}}, \, i=1,2,\cdots, d.   
    \end{align*}
\end{prop}
\begin{proof}
    For any $\theta \in\R^d$,
    Define $h(x;\theta) := e^{\langle\theta,x\rangle}-\langle\theta,x\rangle-1$ and
    \begin{align}
        K^{N}_t(\theta) := \exp{\langle\theta, M^N_t\rangle-\int_0^{t}\sum_{j\in J^N} q^{N}(X^N_{s-},j)h(j;\theta)ds}.
    \end{align}
    It is easy to see that $\left(K^{N}_{t}(\theta)\right)_{t\geq 0}$ is a martingale with mean $1$, since
        \begin{align*}
            & K^{N}_t(\theta)  = \exp{\langle\theta, X^N_t\rangle-\int_0^{t}\sum_{j\in J^N} q^N\left(X^N_{s-},j\right)\left(e^{\langle\theta,j\rangle}-1\right)ds},\\
          = &K^{N}_{0}(\theta)-\int_0^t K^{N}_{s-}(\theta)\sum_{j\in J^N} \left(e^{\langle\theta,j\rangle}-1\right)\upsilon(ds,j) \\
            & + \int_0^t\exp{-\int_0^{s}\sum_{j\in J^N} q^N\left(X^N_{u-},j\right)\left(e^{\langle\theta,j\rangle}-1\right)du}\sum_{j\in J^N}\left(e^{\langle\theta,X^N_{s-}+j\rangle}-e^{\langle\theta,X^N_{s-}\rangle}\right)\mu(ds,j) \\
             = & 1 + \int_0^t K^{N}_{s-}(\theta)\sum_{j\in J^N} \left(e^{\langle\theta,j\rangle}-1\right)(\mu-\upsilon)(ds,j), 
        \end{align*}
        where $K^N(\theta)$ is bounded and $\product{x,y}$ denotes the scalar product.
        
    For any  $x,\theta\in\R^d$,
    \begin{align}\label{ieq:h(x;theta)}
        h(x;\theta)\leq e^{|\product{\theta,x}|}-|\product{\theta,x}|-1\leq \frac{1}{2}|\product{\theta,x}|^2e^{|\product{\theta,x}|}.
    \end{align}
    Conditioning on the event $\{t< \tau^{t_0}(i,\bar{a})\}$ and letting $\theta = ce_i$ for any $c=c(N)>0$, we have 
    \begin{align}\label{assum:expo_drift_exact}
        \int_0^{t\wedge t_0}\sum_{j\in J^N} q^{N}(X^N_{s-},j)h(j;\pm c e_i)ds \leq \frac{1}{2}c^2\exp{c\max_{j\in J^N}\abs{j_i}} \bar{a}, \quad i=1,2,\cdots, d.
    \end{align}
    Now for $t>0$, any $A = A(N)>0$ and $B=B(N)>0$, let \[\tau_M(\theta,B):= \inf\{t:\product{\theta,M^N_t}\geq B\}.\]  
    Since $h(j;\theta)$ is non-negative, 
    \begin{align*}
         & \prob{}{\tau_M(\theta,B) \leq t_0\wedge \tau^{t_0}(i, \bar{a}), \int_0^{\tau_M(\theta,B)\wedge t_0}\sum_{j\in J^N} q^{N}(X^N_{s-},j)h(j;\theta)ds < A} \\
        \leq  & \prob{}{\sup_{s\leq t_0\wedge \tau^{t_0}(i, \bar{a})}\product{\theta,M^N_s}\geq B,\int_0^{\tau_M(\theta,B)\wedge t_0}\sum_{j\in J^N} q^{N}(X^N_{s-},j)h(j;\theta)ds < A} \\
        \leq  &\prob{}{K^{N}_{\tau_M(\theta,B)\wedge t_0}(\theta) \geq e^{B-A}} \\
        \leq & e^{A-B}\expect{}{K^{N}_{\tau_M(\theta,B)\wedge t_0}(\theta)} =  e^{A-B},
    \end{align*}
    where the last equality follows from Doob's optional sampling theorem.

    Let
    \[\theta =  \pm \frac{\delta}{2\bar{a}} e_i.\]
    It follows from \eqref{assum:expo_drift_exact} that on the event $\{t< \tau^{t_0}(i,\bar{a})\}$,
    \[ \int_0^{t\wedge t_0}\sum_{j\in J^N} q^{N}(X^N_{s-},j)h\left(j;\pm \frac{\delta}{2\bar{a}}e_i\right)ds \leq \frac{\delta^2}{4\bar{a}}, \quad i=1,2,\cdots, d.\]
     We then have
    \begin{align*}
        & \prob{}{\sup_{s\leq t_0\wedge  \tau^{t_0}(i,\bar{a})}\abs{\product{e_i,M^N_s}}\geq \delta}\\
        \leq  & \prob{}{\tau_M\left(\frac{\delta}{2\bar{a}} e_i,\frac{\delta^2}{2\bar{a}}\right) \leq t_0\wedge \tau^{t_0}(i,\bar{a})} + \prob{}{\tau_M\left(-\frac{\delta}{2\bar{a}} e_i,\frac{\delta^2}{2\bar{a}}\right) \leq t_0\wedge \tau^{t_0}(i,\bar{a})}\\       
        \leq & 2\exp{\frac{\delta^2}{4\bar{a}}-\frac{\delta^2}{2\bar{a}}} = 2\exp{-\frac{\delta^2}{4\bar{a}}},
    \end{align*}
    and the statement follows.
\end{proof}

\bibliographystyle{abbrvnat}
\bibliography{Library} 
\end{document}

%% file: images/graph_sub.tex
\begin{tikzpicture}[scale=1]
\begin{axis}[
        xmin=0, xmax=1.2,
        ymin=0, ymax=1.2,
        axis lines=middle,
        x label style={at={(axis description cs:1,-0.1)},anchor=north},
         y label style={at={(axis description cs:-0.1,1)},anchor=south},
        xlabel = $\polyorder{I_0}$, ylabel = $\polyorder{R_0}$,
        xtick={0,0.4, 0.45, 0.5, 1},
        xticklabels={0, $a$, $b$, $c$, 1},
        ytick={0,0.6,1},
        yticklabels={0,$1+\polyorder{1-\paralambda}$,1},
      ]  
      \addplot[samples=100, domain=0:1,dash pattern= on 3pt off 1pt, name path=A] {1-x} node[above,pos=1] {$1-\polyorder{I_0}$}; 
      \addplot[samples=50, 
      dash pattern= on 3pt off 1pt,
      domain=0:0.7,name path=B] {x+0.1} node[right,pos=1] {$\polyorder{I_0} -\polyorder{\paragamma}$}; 
      \addplot[samples=50, domain=0:1,name path=Rmax, dotted] {1}; 
      \path[name path=xaxis] (\pgfkeysvalueof{/pgfplots/xmin}, 0) -- (\pgfkeysvalueof{/pgfplots/xmax},0);
       \addplot[fill=lightgray!30] fill between[of=A and xaxis, soft clip={domain=0:0.4}];
      \addplot [ultra thick,mark=none,solid] coordinates {(0.4, 0) (0.4, 0.6)};
      \addplot [mark=none,dash pattern= on 3pt off 1pt, name path = C] coordinates {(0, 0.6) (1, 0.6)};
      
       \addplot[fill=lightgray!75] fill between[of=C and xaxis, soft clip={domain=0.4:0.5}];
       \addplot [ mark=none, 
       dash pattern= on 3pt off 1pt,
       ] coordinates {(0.45, 0) (0.45, 0.55)};
       \addplot [mark=none, 
       dash pattern= on 3pt off 1pt,
       ] coordinates {(0.5, 0) (0.5, 0.6)};
      \node at (0.2,0.5) {\small$1.1$};
      \node at (0.4,0.8) {\small$1.2$};
      \draw[->] (0.4,0.75) -- (0.4,0.65);
      \node at (0.45,0.3) {\small$1.3$};
      \addplot [ultra thick,mark=none,solid] coordinates {(0, 1) (0.95, 1)};
       \node at (0.5,1.15) {\small$2.2$};
       \draw[->] (0.45,1.08) -- (0.4,1.02);
        \node at (0.2,1.15) {\small$2.1$};
      \draw[->] (0.15,1.1) -- (0.02,1.02);
      \node at (0,1)[circle,fill,inner sep=1.5pt]{};
   \end{axis}
\end{tikzpicture}
\hfill
\begin{tikzpicture}[scale=1]
\begin{axis}[
        xmin=0, xmax=1.2,
        ymin=0, ymax=1.2,
        axis lines=middle,
        x label style={at={(axis description cs:1,-0.1)},anchor=north},
         y label style={at={(axis description cs:-0.1,1)},anchor=south},
        xlabel = $\polyorder{I_0}$, ylabel = $\polyorder{R_0}$,
        xtick={0,0.2, 0.3, 0.4, 1},
        xticklabels={0, $c$, $b$, $a$,1},
        ytick={0,0.6,1},
        yticklabels={0,$1+\polyorder{1-\paralambda}$,1},
      ]  
      \addplot[samples=50, domain=0:1,name path=Rmax, dotted] {1}; 
      \path[name path=xaxis] (\pgfkeysvalueof{/pgfplots/xmin}, 0) -- (\pgfkeysvalueof{/pgfplots/xmax},0);
    
      \addplot[samples=100, domain=0:1,dash pattern= on 3pt off 1pt, name path=A] {1-x} node[above,pos=1] {$1-\polyorder{I_0}$}; 
     
      \addplot[samples=50, 
      dash pattern= on 3pt off 1pt,
      domain=0:0.6,name path=B] {x+0.4} node[right,pos=0.8] {$\polyorder{I_0} -\polyorder{\paragamma}$}; 
      
       \addplot[fill=lightgray!30] fill between[of=A and xaxis, soft clip={domain=0:0.3}];
      
      \addplot [mark=none, dash pattern= on 3pt off 1pt, name path = C] coordinates {(0, 0.6) (1, 0.6)};
      
       \addplot [mark=none, dash pattern= on 3pt off 1pt] coordinates {(0.2, 0) (0.2, 0.6)};
       \addplot [mark=none, 
       ] coordinates {(0.3, 0) (0.3, 0.7)};
       \addplot [mark=none, 
       dash pattern= on 3pt off 1pt,
       ] coordinates {(0.4, 0) (0.4, 0.6)};
      \node at (0.15,0.3) {\small$1.1$};
      \addplot [ultra thick,mark=none,solid] coordinates {(0, 1) (0.95, 1)};
       \node at (0.45,1.15) {\small$2.2$};
       \draw[->] (0.45,1.08) -- (0.4,1.02);
        \node at (0.2,1.15) {\small$2.1$};
      \draw[->] (0.15,1.1) -- (0.02,1.02);
      \node at (0,1)[circle,fill,inner sep=1.5pt]{};
   \end{axis}
\end{tikzpicture}